\documentclass{amsart}
\usepackage{amsmath,amsfonts,amssymb,amsthm,amscd, mathrsfs}
\usepackage{graphics}
\usepackage{ae, aecompl, aeguill}
\usepackage{url}
\usepackage{hyperref}

\newtheorem{thm}[equation]{Theorem}
\newtheorem{lem}[equation]{Lemma}
\newtheorem{pro}[equation]{Proposition}
\newtheorem{cor}[equation]{Corollary}
\theoremstyle{definition}
\newtheorem{defi}[equation]{Definition}
\newtheorem{nota}[equation]{Notation}
\newtheorem{rem}[equation]{Remark}

\DeclareMathOperator\Aut{Aut}

\DeclareMathOperator\GL{GL}
\DeclareMathOperator\dep{dp}
\DeclareMathOperator\gr{gp}
\DeclareMathOperator\Image{Im}
\DeclareMathOperator\Id{Id}

\newcommand{\N}{\mathbb{N}}

\newcommand{\R}{\mathbb{R}}

\newcommand{\K}{\mathbb{K}}
\newcommand{\EL}{\mathscr{L}}
\newcommand{\gras}[1]{\textrm{\boldmath $#1$}}
\newcommand{\ie}{{\em i.e. }}

\newcommand{\cf}{{\em cf. }}
\def\G{\Gamma }
\def\D{\Delta }

\def\S{\Sigma }
\def\a{\alpha }
\def\b{\beta }

\def\s{\sigma }

\def\p{\varphi }
\def\m{\mathfrak }
\def\bt{\otimes }

\begin{document}
\author{Anatole Castella}
\title[Twisted Lawrence-Krammer representations]{Twisted Lawrence-Krammer representations}
\date{\today}

\maketitle

\begin{abstract}
Lawrence-Krammer representations are an important family of linear representations of Artin-Tits groups of small type, which are known, under some assumptions on the parameters, to be faithful when the type is spherical (or more generally when they are restricted to the Artin-Tits monoid) and irreducible when the type is connected. 

Here, we investigate an analogue of these representations --- introduced by Digne in the spherical cases --- for every Artin-Tits monoid that appears as the submonoid of fixed points of an Artin-Tits monoid of small type under a group of graph automorphisms, and for the corresponding Artin-Tits group. 

Under the same assumptions on the parameters as in the small type cases, we first show that these so-called ``twisted Lawrence-Krammer representations'' are faithful, and we then prove, by computing their formulas when the group of graph automorphisms is of order two or three, their irreducibility in all the spherical and connected cases but one. 
\end{abstract}

\section*{Introduction}

Lawrence-Krammer representations (shortened to LK-representations in what follows) are a family of linear representations of Artin-Tits monoids and groups of small type, introduced by Lawrence \cite{L} and Krammer \cite{K} for the braid groups, and intensively studied since then (see the references listed below). 

Under some assumptions on the defining ring and on the parameters (see condition ($\bigstar$) of theorem \ref{main thm 1} below), they are known to be faithful for the monoids \cite{K2,B,CW,Di,P,H,Ca} which ensure their faithfulness for the groups when the type is spherical, irreducible when the type is connected \cite{Zin,Ma,CGW,Ca}, and have proved to be useful in the study of several other properties of Artin-Tits groups \cite{Ma} and related objects \cite{CGW,Ma2}. 

\medskip

It is therefore an interesting question to ask if there exists some analogous linear representations for Artin-Tits monoids and groups of non-small type.

\medskip

In \cite[end of section 3]{Di}, Digne defines such objects for Artin-Tits monoids and groups of spherical type $B_n$, $F_4$ and $G_2$, using the fact that they appear as the submonoids and subgroups of fixed points of Artin-Tits monoids and groups of small and spherical type (namely of type $A_{2n-1}$, $E_6$ and $D_4$ respectively) under the action of a graph automorphism. 

These new linear representations --- that I call \emph{``twisted'' LK-representations} --- share the same good combinatorial properties with respect to the associated root systems as in the small type cases, and Digne proves their faithfulness for his choice of parameters in \cite[Cor. 3.11]{Di}.

\medskip

The aim of this article is to generalize this construction to every Artin-Tits monoid that appears as the submonoid of fixed elements of an Artin-Tits monoid of small type under a group of graph automorphisms, and to the associated Artin-Tits group. 

Under the same condition ($\bigstar$) on the defining ring and on the parameters as in the small type cases, we prove their faithfulness for the monoids, which again ensure their faithfulness for the groups when the type is spherical, and we prove their irreducibility for all the connected and spherical types but one. 

\medskip

The paper is organized as follows. 

We recall the basic needed notions on Coxeter groups, Artin-tits monoids an groups, standard root systems and graph automorphisms in section \ref{Preliminaries}. 

We recall the definition and the main properties of LK-representations in the small type cases, as stated in \cite{Ca}, in subsection \ref{Lawrence-Krammer representations}, and turn to the definition of the twisted LK-representations in subsection \ref{Definition twisted}. We prove our faithfulness result in subsection \ref{Twisted faithfulness criterion} (theorem \ref{twisted criterion}).

We explicit the formulas of these twisted LK-representations when the group of graph automorphisms is of order two or three in section \ref{Case two or three}.

We use these formulas in section \ref{The spherical case} to study the spherical cases. We prove our irreducibility result in subsection \ref{Irreducibility} (theorem \ref{thm irred}), and conclude in subsection \ref{Non-equivalence} by giving a sufficient condition on the parameters for non-equivalence in the family of twisted LK-representations of a given type.

\section{Coxeter matrices and related objects}\label{Preliminaries}

A \emph{Coxeter matrix} is a matrix $\G = (m_{i,j})_{i,j\in I}$ over an arbitrary set $I$ with $m_{i,j} = m_{j,i} \in \N_{\geqslant 1}\cup\{\infty\}$, and $m_{i,j} = 1 \Leftrightarrow i = j$ for all $i,\,j\in I$. 

As usual, we encode the data of $\G$ by its \emph{Coxeter graph}, \ie the graph with vertex set $I$, an edge between the vertices $i$ and $j$ if $m_{i,j} \geqslant 3$, and a label $m_{i,j}$ on that edge when $m_{i,j} \geqslant 4$. 

We say that a Coxeter matrix $\G$ is \emph{connected} when so is its Coxeter graph.

\medskip

In this paper, we will always assume that $I$ is finite. This condition could be removed at the cost of some refinements in certain statements below (see \cite[Ch.~11]{Ca2} for some of them), which are left to the reader.

\subsection{Coxeter groups and Artin-Tits monoids and groups}\label{Main notations and definitions}\mbox{}\medskip

To a Coxeter matrix $\G = (m_{i,j})_{i,j\in I}$, we associate the \emph{Coxeter group} $W = W_\G$, the \emph{Artin-Tits group} $B = B_\G$ and the \emph{Artin-Tits monoid} $B^+ = B^+_\G$, given by the following presentations : 
\begin{eqnarray*}
W & = & \langle\ s_i,\, i\in I\mid \underbrace{s_i s_j s_i \cdots}_{m_{i,j}\, \text{terms}} = \underbrace{s_j s_i s_j \cdots}_{m_{i,j}\, \text{terms}} \ \text{   if   }\  m_{i,j}\neq \infty, \text{   and   } s_i^2 = 1 \ \rangle,\\
B & = & \langle\ \gras s_i,\, i\in I\mid \underbrace{\gras s_i\gras s_j\gras s_i \cdots}_{m_{i,j}\, \text{terms}} = \underbrace{\gras s_j\gras s_i\gras s_j \cdots}_{m_{i,j}\, \text{terms}} \  \text{   if   }\ m_{i,j}\neq \infty\ \rangle,\\
B^+ & = & \langle\ \gras s_i,\, i\in I\mid \underbrace{\gras s_i\gras s_j\gras s_i \cdots}_{m_{i,j}\, \text{terms}} = \underbrace{\gras s_j\gras s_i\gras s_j \cdots}_{m_{i,j}\, \text{terms}} \  \text{   if   }\ m_{i,j}\neq \infty\ \rangle^+.
\end{eqnarray*}

We denote by $\ell$ the \emph{length} function on $W$ and on $B^+$ relatively to their generating sets $\{s_i, \, i\in I\}$ and $\{\gras s_i, \, i\in I\}$ respectively. 

\medskip

We denote by $\preccurlyeq$ the \emph{left divisibility} in the monoid $B^+$, \ie for $g, \,h \in B^+$, we write $h \preccurlyeq g$ if there exists $h'\in B^+$ such that $g = hh'$. This leads to the natural notions of \emph{left gcd}'s and \emph{right lcm}'s in $B^+$. 

Let $J$ be a subset of $I$. We denote by $\G_J$ the submatrix $(m_{i,j})_{i,j\in J}$ of $\G$ and by $W_J$ the subgroup $\langle s_j,\, j\in J\rangle$ of $W$. We say that $J$ and $\G_J$ are \emph{spherical} if $W_J$ is finite, or equivalently if the elements $\gras s_j$, $j \in J$, have a common right multiple in $B^+$ (see \cite[Thm. 5.6]{BS}). In that case, there exists a unique element $r_J$ of maximal length in $W_J$, and the elements $\gras s_j$ ($j \in J$) have a unique right lcm in $B^+$ that we denote by $\D_J$ (see \cite[Props. 4.1 and 5.7]{BS}). 

\subsection{Standard root system}\label{Standard root system}\mbox{}\medskip

We refer the reader to \cite{Deo} for the basic notions on standard root systems. 

We denote by $\Phi = \Phi_\G = \{w(\a_i) \mid w\in W, \, i\in I\}$ the \emph{standard root system} associated with $\G$ in the real vector space $E = \oplus_{i\in I}\R\a_i$, where the action of $W$ on $E$ is defined, for a generator $s_i$ and a basis element $\a_j$, by $s_i(\a_j) = \a_j+2\cos\big(\frac{\pi}{m_{i,j}}\big)\a_i$. 

It is known that $\Phi = \Phi^+ \sqcup \Phi^-$, where $\Phi^+ = \Phi \bigcap \left(\oplus_{i\in I}\R^+\a_i\right)$ and $\Phi^- = -\Phi^+$. 

\medskip

As in \cite{Ca}, we will represent $\Phi^+$ --- and any of its subsets --- as a graded graph, where two elements $\a$ and $\b$ of $\Phi^+$ are linked by an edge labeled $i$ if $\a = s_i(\b)$, and where the grading is by the \emph{depth} function on $\Phi^+$ defined, for $\a\in \Phi^+$, by $$\dep(\a) = \min\{\ell(w) \mid w \in W\, \textrm{ such that }\, w(\a) \in \Phi^-\}.$$ By convention in those graded graphs, we choose to place roots of great depth \emph{above} roots of small depth, so edges like the following ones will all mean that $\a = s_i(\b)$ and $\dep(\b) > \dep(\a)$ : 
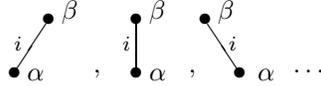
\begin{center}
\begin{picture}(30,25)
 \put(0,0){\circle*{4}}\put(5,-3){$\a$}
 \put(13,20){\circle*{4}}\put(18,20){$\b$}
 \put(0,0){\line(2,3){12}}\put(0,8){\footnotesize $i$}
\end{picture},
\begin{picture}(30,25)
 \put(10,0){\circle*{4}}\put(15,-3){$\a$}
 \put(10,20){\circle*{4}}\put(15,20){$\b$}
 \put(10,0){\line(0,1){20}}\put(5,8){\footnotesize $i$}
\end{picture},
\begin{picture}(30,25)
 \put(13,0){\circle*{4}}\put(20,-3){$\a$}
 \put(0,20){\circle*{4}}\put(5,20){$\b$}
 \put(13,0){\line(-2,3){12}}\put(9,8){\footnotesize $i$}
\end{picture}
$\ldots$
\end{center}

\medskip

In particular in section \ref{Case two or three}, we will focus on some particular subsets of $\Phi^+$, closure of subsets of $\Phi^+$ under the action of some subgroup $W_{J}$ of $W$, that we call, after H\'ee, \emph{meshes} :

\begin{nota}For $J\subseteq I$, we call \emph{$J$-mesh} of a subset $\Psi$ of $\Phi^+$ the set $M_{J}(\Psi) = W_{J}(\Psi)\bigcap \Phi^+ =  \{w(\a) \mid w\in W_{J},\, \a\in \Psi\, \textrm{ such that }\, w(\a)\in \Phi^+\}$. 
\end{nota}

\subsection{Graph automorphisms}\label{Graph automorphisms}\mbox{}\medskip

We call {\em graph automorphism} of a Coxeter matrix $\G = (m_{i,j})_{i,j\in I}$ every permutation $\s$ of $I$ such that $m_{\s(i),\s(j)} = m_{i,j}$ for all $i,\,j \in I$, and we denote by $\Aut(\G)$ the group they constitute.

\medskip

Any graph automorphism $\s$ of $\G$ acts by an automorphism on $W$ (resp. on $B$ and $B^+$) by permuting the generating set $\{s_i \mid i\in I\}$ (resp. $\{\gras s_i \mid i\in I\}$). 

If $\S$ is a subgroup of $\Aut(\G)$, we denote by $W^\S$ (resp. $B^\S$, resp. $(B^+)^\S$) the subgroup of $W$ (resp. subgroup of $B$, resp. submonoid of $B^+$) of fixed points under the action of the elements of $\S$. It is known that $W^\S$ (resp. $(B^+)^\S$) is a Coxeter group (resp. Artin-Tits monoid), and that the analogue holds for $B^\S$ when $\G$ is spherical, or more generally of FC-type (see \cite{H1,Mu} for the Coxeter case, \cite{Mi,C2,C2',Ca1} for the Artin-Tits case). 

More precisely, if we denote by $I^\S$ the set of spherical orbits of $I$ under $\S$, then $W^\S$ (resp. $(B^+)^\S$, and $B^\S$ when $\G$ is of FC-type) has a Coxeter (resp. Artin-Tits) presentation associated with a Coxeter matrix $\G^\S = (m^\S_{J,K})_{J,K \in I^\S}$ easily computable from $\G$, where the generators are the elements $r_J$ (resp. $\D_J$) for $J$ running through $I^\S$ (see for example \cite{C2,C2',Ca1}).
\medskip

Similarly, any graph automorphism $\s$ of $\G$ acts by a linear automorphism on the vector space $E = \oplus_{i\in I}\R\a_i$ by permuting the basis $(\a_i)_{i\in I}$. This action stabilizes $\Phi$ and $\Phi^+$, and the induced action on those sets is $w(\a_i) \mapsto (\s(w))(\a_{\s(i)})$. 

It is easily seen on the definition that the action of $\Aut(\G)$ on $\Phi^+$ thus defined respects the depth function on $\Phi^+$. In particular, we can define the \emph{depth} of an orbit $\Theta$ of $\Phi^+$ under $\S$ as the depth of any element $\a \in \Theta$.

\section{Twisted Lawrence-Krammer representations}\label{Action of graph automorphisms}

From now on, we fix a Coxeter matrix $\G = (m_{i,j})_{i,j\in I}$ of \emph{small type}, \ie with $m_{i,j} \in \{2,3\}$ for all $i,j\in I$ with $i\neq j$.

\subsection{Lawrence-Krammer representations - the small type case}\label{Lawrence-Krammer representations}\mbox{}\medskip

%

Let $R$ be a (unitary) commutative ring and $V$  be a free $R$-module with basis $(e_\a)_{\a \in \Phi^+}$. We denote by $V^\star$ the dual of $V$ and by $R^\times$  the group of units of $R$. 

For $f \in V^\star$ and $e \in V$, we denote by $f\bt e$ the endomorphism of $V$ defined by $v \mapsto f(v)e$.

\begin{defi}[{\cite[Def. 7]{Ca}}]\label{def psii}For ${(a,b,c,d) \in R}^4$ and $i \in I$, we denote by $\p_i = \p_{i,(a,b,c,d)}$ the endomorphism of $V$ given on the basis $(e_\a)_{\a\in \Phi^+}$ by 
\begin{center}
$\begin{cases}
 \p_i(e_{\a}) = 0 & \text{  if  }\ \ \a = \a_i, \\
\p_i(e_\a) = d e_\a & \text{  if  
\begin{picture}(36,12)(0,-2)
\put(20,0){\circle*{4}}
\put(25,0){\circle{10}}
\put(32,-3){\footnotesize $i$}
\put(10,2){$\a$}
\end{picture}},\\
\begin{cases}\p_i(e_\b) = b e_\a \\ \p_i(e_\a) = a e_\a + c e_\b \end{cases} & \text{  if  
\begin{picture}(32,20)(0,8)
\put(20,0){\circle*{4}}
\put(20,20){\circle*{4}}
\put(20,0){\line(0,1){20}}
\put(22,8){\footnotesize $i$}
\put(10,0){$\a$}
\put(10,18){$\b$}
\end{picture} in   } \, \Phi^+.
\end{cases}$
\end{center}

For a linear form $f_i \in V^\star$, we define the \emph{Lawrence-Krammer map} --- or the \emph{LK-map} for short ---  associated with ${(a, b, c, d)}$ and $f_i$ to be the endomorphism $\psi_i = \psi_{i,(a,b,c,d),f_i}$ of $V$ given by $\psi_i = \p_i+f_i\bt e_{\a_i}$.
\end{defi}

\begin{pro}[{\cite[Lem. 9 and Prop. 12]{Ca}}]\label{prop LK-trucs}
Assume that $d^2 - ad - bc = 0$ and that the family $(f_i)_{i\in I} \in (V^\star)^I$ satisfies the following properties :
\begin{enumerate}
\item for $i,\,j \in I$ with $i\neq j$, $f_i(e_{\a_j})=0$,
\item for $i,\, j \in I$ with $m_{i,j} = 2$, $f_i\p_j = df_i$,
\item for $i,\, j \in I$ with $m_{i,j} = 3$, $f_i\p_j = f_j\p_i$.
\end{enumerate}
Then $\gras s_i \mapsto \psi_i$ defines a linear representation $\psi : B^+ \to \EL(V)$. Moreover if $b$, $c$, $d$ and $f_i(e_{\a_i})$, $i\in I$, belong to $R^\times$, then the image of $\psi$ is included in $\GL(V)$ and hence $\psi$ induces a linear representation $\psi_{\gr} : B \to \GL(V)$.
\end{pro}

\begin{defi}[{\cite[Def. 11]{Ca}}]\label{LK-trucs}We call \emph{LK-family (relatively to $(a,b,c,d)$)} any family $(f_i)_{i\in I} \in (V^\star)^I$ satisfying conditions (i), (ii) and (iii) of proposition \ref{prop LK-trucs} above, and we call \emph{LK-representation (relatively to $(a,b,c,d)$ and $(f_i)_{i\in I}$)} the induced representation $\psi$ and, when appropriate, $\psi_{\gr}$.
\end{defi}

\begin{rem}
When $b$ is invertible in $R$, an LK-family $(f_i)_{i\in I}$ necessarily satisfies $f_i(e_{\a_i}) = f_j(e_{\a_j})$ for every $i$, $j\in I$ with $m_{i,j} = 3$ (see \cite[Prop. 34]{Ca}). In particular when $\G$ is connected, we thus get that the elements $f_i(e_{\a_i})$, $i\in I$, are all equal. For sake of brevity when this is the case, we will simply denote by $f$ the common value of the $f_i(e_{\a_i})$, $i\in I$.

Moreover when $\G$ is connected and spherical, one can show that the LK-family $(f_i)_{i\in I}$ is entirely determined by this common value $f$ (see \cite[Section 3.2]{Ca}). 
\end{rem}

The first main properties of these representations can be stated as follows :

\begin{thm}[{\cite[Thm. A]{Ca}}]\label{main thm 1}Let $\psi : B^+ \to \EL(V)$ be an LK-representation over $R$ associated with parameters $(a,b,c,d)\in R^4$ and $(f_i)_{i\in I} \in (V^\star)^I$. Assume that the following condition holds : 
\begin{itemize}
\item[($\bigstar$)] $R$ is an integral domain, $f_i(e_{\a_i})\neq 0$ for all $i \in I$, and there exists a totally ordered integral domain $R_0$ and a ring homomorphism $\rho : R \to R_0$, such that $\rho(t) > 0$ for $t \in \{{a,b,c,d}\}$, and $\Image(f_i) \subseteq \ker(\rho)$ for all $i \in I$.
\end{itemize}
By extension of the scalars from $R$ to its field of fractions $\K$, consider $\psi$ as acting on the $\K$-vector space $V_\K = \K\otimes_RV$. Then $\psi$ has its image included in $\GL(V_\K)$, hence induces an LK-representation $\psi_{\gr} : B \to \GL(V_\K)$, and the following holds.
\begin{enumerate}
\item The LK-representation $\psi$ is faithful, and so is $\psi_{gr}$ if $\G$ is spherical.
\item If $\G$ is connected, the LK-representations $\psi$ and $\psi_{\gr}$ are irreducible on $V_\K$.
\item Assume that $\G$ has at least one edge and that $\psi'$ is an LK-representation of $B^+$ associated with parameters $(a',b',c',d') \in R^4$ and $(f_i')_{i\in I}\in (V^\star)^I$ that satisfy the conditions of the preamble for the fixed $\rho$. Then if $a' \neq a$, or if $d' \neq d$ or if $f_i'(e_{\a_i})\neq f_i(e_{\a_i})$ for some $i \in I$, the LK-representations $\psi$ and $\psi'$ are not equivalent on $V_\K$.
\end{enumerate}
\end{thm}

\subsection{Definition of the twisted Lawrence-Krammer representations}\label{Definition twisted}\mbox{}\medskip

Let us fix a subgroup $\S$ of $\Aut(\G)$ and an LK-representation $\psi : B^+\to \EL(V)$, $g \mapsto \psi_g$, of $B^+$ associated with an LK-family $(f_i)_{i\in I} \in (V^\star)^I$. 

Recall that $\S$ naturally acts on $B^+$, and that the submonoid $(B^+)^\S$ of fixed points of $B^+$ under $\S$ is the Artin-Tits monoid $B^+_{\G^\S}$ associated with a certain Coxeter matrix $\G^\S$ (see subsection \ref{Graph automorphisms}). 

Moreover, $\S$ acts on $\Phi^+$ and this action induces an action of $\S$ on $V$ by permutation of the basis $(e_\a)_{\a \in \Phi^+}$, which we denote by $\S \to \GL(V)$, $\s \mapsto \s_V$. We denote by $V^\S$ the submodule of fixed points of $V$ under the action of $\S$.  

\medskip

Let us denote by $\p : B^+ \to \mathscr{L}(V)$, $g \mapsto \p_g$, the LK-representation of $B^+$ associated with the trivial LK-family (\ie where $f_i$ is the zero form for all $i\in I$).

\begin{lem}\label{stabilisation phi}For all $(g,\s) \in B^+ \times \Aut(\G)$, we get $\s_V\p_g = \p_{\s(g)}\s_V$. In particular, for every $g \in (B^+)^\S$, $\p_g$ stabilizes $V^\S$ and hence $\p$ induces a linear representation $$\p^\S : (B^+)^\S \to \mathscr{L}(V^\S), \ \ g \mapsto \p^\S_g = \p_g|_{V^\S}.$$ 
\end{lem}
\proof The action of $\Aut(\G)$ on $\Phi^+$ respects the depth, and this clearly implies that $\s_V(\p_{i}(e_\a)) = \p_{\s(i)}(e_{\s(\a)})$ in view of the formulas of definition \ref{def psii}. The result follows by linearity and induction on $\ell(g)$.\qed

\begin{pro}\label{stabilisation}Assume that $f_{i} = f_{\s(i)}\s_V$ in $V^\star$, for every $(i,\s) \in I \times \S$. Then for every $(g,\s) \in B^+\times \S$, we get $\s_V\psi_g = \psi_{\s(g)}\s_V$. In particular, for every $g \in (B^+)^\S$, $\psi_g$ stabilizes $V^\S$ and hence $\psi$ induces a linear representation $$\psi^\S : (B^+)^\S \to \mathscr{L}(V^\S), \ \ g \mapsto \psi^\S_g = \psi_g|_{V^\S}.$$

Moreover if the images of $\psi$ are invertible, then so are the images of $\psi^\S$. 
\end{pro}
\proof We have $\s_V\psi_i = \s_V\p_i + \s_V(f_i\bt e_{\a_i}) = \s_V\p_i + f_i\bt e_{\a_{\s(i)}}$ and $\psi_{\s(i)}\s_V = \p_{\s(i)}\s_V + (f_{\s(i)}\bt e_{\a_{\s(i)}})\s_V = \p_{\s(i)}\s_V + (f_{\s(i)}\s_V)\bt e_{\a_{\s(i)}}$ (thanks to the formulas of \cite[Rem. 1]{Ca}), whence $\s_V\psi_i = \psi_{\s(i)}\s_V$ by the previous lemma and assumption on $f_i$ and $f_{\s(i)}$. The first point follows by induction on $\ell(g)$. Moreover if the images of $\psi$ are invertible, then the equality $\s_V\psi_g\s_V^{-1} = \psi_{\s(g)}$ implies $\s_V\psi_g^{-1}\s_V^{-1} = \psi_{\s(g)}^{-1}$, and hence $\psi^{-1}_g$ stabilizes $V^\S$ for every $g \in (B^+)^\S$. This gives the result. \qed

\begin{defi}[twisted LK-representations]
Under the assumption of the previous proposition, we call \emph{twisted LK-representation} the linear representation $\psi^\S : (B^+)^\S \to \EL(V^\S)$ of the Artin-Tits monoid $(B^+)^\S = B^+_{\G^\S}$, and, when appropriate, the induced representation $\psi^\S_{\gr}: B_{\G^\S} \to \GL(V^\S)$ of the Artin-Tits group $B_{\G^\S}$.
\end{defi}

The assumption $f_i = f_{\s(i)}\s_V$ for every $(i,\s) \in I \times \S$ is equivalent to $f_{i}(e_\a) = f_{\s(i)}(e_{\s(\a)})$ for every $(i,\a,\s) \in I\times \Phi^+ \times \S$. It is not always satisfied : for example if $i$ and $\s(i)$ are not in the same connected component of $\G$, then $f_i(e_{\a_i})$ and $f_{\s(i)}(e_{\a_{\s(i)}})$ can be chosen to be distinct (see \cite[Section 3.1]{Ca}). I do not know if this assumption is always satisfied when $\G$ is connected, but we have the following partial result : 

\begin{pro}\label{cas favorables}Let $(f_i)_{i\in I}$ be an LK-family. Assume that $b$, $c$, $d$ are invertible in $R$ and that we are in one of the following cases : 
\begin{enumerate}
 \item $\G$ is spherical and connected (\ie of type $ADE$), or
 \item $\G$ is affine (\ie of type $\tilde A\tilde D\tilde E$), or
 \item $\G$ has \emph{no triangle} and $(f_i)_{i\in I}$ is the LK-family of Paris, as in \cite[Def. 42]{Ca}.
 \end{enumerate}

Then $f_i = f_{\s(i)}\s_V$ for every $(i,\s) \in I \times \Aut(\G)$.
\end{pro}
\proof The condition $f_{i}(e_\a) = f_{\s(i)}(e_{\s(\a)})$ for every $(i,\a,\s) \in I\times \Phi^+\times \Aut(\G)$, for the three situations (note that the first one is a consequence of the third one), is easy to see by induction on $\dep(\a)$, using the inductive construction of the $f_{i}(e_\a)$, $(i,\a)\in I\times \Phi^+$, and the independence results at the inductive steps, of \cite[Sections 3.2, 3.3 and 3.4]{Ca} respectively, and using the fact that the action of $\Aut(\G)$ on $\Phi^+$ respects the depth. \qed

\subsection{Twisted faithfulness criterion}\label{Twisted faithfulness criterion}\mbox{}\medskip

The aim of this subsection is to prove that the classical faithfulness criterion of \cite{H}, as stated in theorem \ref{main thm 1} above, also works in the twisted cases. 

\medskip

For $g \in B^+$, we set $I(g) = \{i\in I \mid \gras s_i \preccurlyeq g\}$. 

\begin{lem}\label{reduction2}Let $\psi : (B^+)^\S \to M$ be a monoid homomorphism where $M$ is left cancellative. If $\psi$ satisfies $\psi(g) = \psi(g') \Rightarrow I(g)= I(g')$ for all $g,\,g' \in (B^+)^\S$, then $\psi$ is injective.
\end{lem}
\proof Let $g$, $g' \in (B^+)^\S$ be such that $\psi(g) = \psi(g')$. We prove by induction on $\ell(g)$ that $g = g'$. If $\ell(g) = 0$, \ie if $g = 1$, then $I(g)= I(g') = \emptyset$, hence $g' = 1$ and we are done. If $\ell(g) > 0$, fix $i \in I(g) = I(g')$. Since the action of $\S$ on $B^+$ respects the divisibility and since $g$ is fixed by $\S$, the orbit $J$ of $i$ under $\S$ is included in $I(g) = I(g')$, but then $J$ is spherical and there exist $g_1$, $g'_1\in B^+$ such that $g = \D_Jg_1$ and $g' = \D_Jg'_1$. Since the elements $g$, $g'$ and $\D_J$ are fixed by $\S$ (recall that $\D_J$ is a generator of $(B^+)^\S$), so are $g_1$ and $g'_1$ by cancellation in $B^+$. We thus get $\psi(\D_J)\psi(g_1) = \psi(\D_J)\psi(g'_1)$ in $M$, whence $\psi(g_1) = \psi(g'_1)$ by cancellation in $M$. We therefore get $g_1 = g'_1$ by induction and finally $g = g'$. \qed

\begin{nota}We denote by $\Phi^+/\S$ the set of orbits of $\Phi^+$ under $\S$ and, for every $\Theta \in \Phi^+/\S$, we set $e_\Theta = \sum_{\a\in \Theta}e_\a$. The family $(e_\Theta)_{\Theta \in \Phi^+/\S}$ is a basis of $V^\S$. 
\end{nota}

\begin{thm}\label{twisted criterion}Assume that $f_{i} = f_{\s(i)}\s_V$ for every $(i,\s) \in I \times \S$, so that the twisted LK-representation $\psi^\S : (B^+)^\S \to \EL(V^\S)$ is defined (see proposition \ref{stabilisation}), and assume that condition \emph{($\bigstar$)} holds : 
\begin{itemize}
\item[($\bigstar$)] $R$ is an integral domain, $f_i(e_{\a_i})\neq 0$ for all $i \in I$, and there exists a totally ordered integral domain $R_0$ and a ring homomorphism $\rho : R \to R_0$, such that $\rho(t) > 0$ for $t \in \{{a,b,c,d}\}$, and $\Image(f_i) \subseteq \ker(\rho)$ for all $i \in I$.
\end{itemize}
By extension of the scalars from $R$ to its field of fractions $\K$, consider $\psi^\S$ as acting on the $\K$-vector space $V^\S_\K = \K\otimes_RV^\S$, so that $\psi^\S$ has its image included in $\GL(V^\S_\K)$, and hence induces a twisted LK-representation $\psi^\S_{\gr} : B_{\G^\S} \to \GL(V^\S_\K)$. Then the twisted LK-representation $\psi^\S$ is faithful, and so is $\psi^\S_{\gr}$ if $\G^\S$ is spherical.
\end{thm}
\proof Under condition ($\bigstar$), the elements $b$, $c$, $d$ and $f_i(e_{\a_i})$, $i \in I$, are non-zero in $R$, so become units of the field of fractions $\K$ of $R$. By propositions \ref{prop LK-trucs} and \ref{stabilisation}, $\Image(\psi^\S)$ is then included in $\GL(V^\S_\K)$ --- hence is cancellative --- and $\psi^\S_{\gr} : B_{\G^\S} \to \GL(V^\S_\K)$ is defined. Moreover when $\G^\S$ is spherical, the faithfulness of $\psi^\S$ implies the one of $\psi^\S_{\gr}$ by \cite[Lem. 6]{Ca}. 

In order to prove the theorem, it then suffices to see that $\psi^\S$ satisfies the assumption of lemma \ref{reduction2}. So let $g,\, g' \in (B^+)^\S$ be such that $\psi^\S_g = \psi^\S_{g'}$ and let us show that $I(g) = I(g')$. We need for that some notations of \cite[Section 2.2]{Ca}. 

If we denote by $V_0$ the free $R_0$-module with basis $(e_\a)_{\a \in \Phi^+}$, then the morphism $\rho : R \to R_0$ induces a natural monoid homomorphism $\tilde\rho : \EL(V)\to \EL(V_0)$, $\p \mapsto \overline{\p}$. By ($\bigstar$), the homomorphism $\tilde\rho$ sends $\Image(\psi)$ into $\EL^+(V_0)$, the submonoid of $\EL(V_0)$ composed of the endomorphisms of $V_0$ whose matrix in the basis $(e_\a)_{\a\in \Phi^+}$ has non-negative coefficients. Now for $h\in (B^+)^\S$ and $\b \in \Phi^+$, let us denote by $\m R_h(\b)$ the support of $\overline{\psi_h}(e_\b)$ in the basis $(e_\a)_{\a\in \Phi^+}$, and for $\Psi \subseteq \Phi^+$, let us set $\m R_h(\Psi) = \bigcup_{\b\in \Psi} \m R_h(\b)$. Since the coefficients of the matrix of $\overline{\psi_{h}}$ in the basis $(e_\a)_{\a \in \Phi^+}$ of $V_0$ are non-negative, the set $\m R_h(\Psi)$, when $\Psi$ is finite, is precisely the support of $\overline{\psi_h}(\sum_{\b\in \Psi}e_\b)$ in the basis $(e_\a)_{\a \in \Phi^+}$. 

In order to show that $I(g) = I(g')$, it suffices to prove, in view of \cite[Prop.~2]{H} (see \cite[Lem.~20]{Ca}), that $\m R_g(\Phi^+) = \m R_{g'}(\Phi^+)$. But since $\psi_{g}$ and $\psi_{g'}$ coincide on $V^\S$, we get in particular that $\psi_g(e_\Theta) = \psi_{g'}(e_\Theta)$, and hence $\m R_g(\Theta) = \m R_{g'}(\Theta)$ for every $\Theta \in \Phi^+/\S$. And finally since $\Phi^+ = \bigcup_{\Theta \in \Phi^+/\S}\Theta$, we get $\m R_g(\Phi^+) = \bigcup_{\Theta \in \Phi^+/\S}\m R_g(\Theta) = \bigcup_{\Theta \in \Phi^+/\S}\m R_{g'}(\Theta) = \m R_{g'}(\Phi^+)$, whence the result. \qed



\section{Case of an automorphism of order two or three}\label{Case two or three}

Let $\G = (m_{i,j})_{i,j\in I}$ be a Coxeter matrix of small type and fix $\S = \langle\s \rangle \leqslant \Aut(\G)$ of order two or three.

Recall that the submonoid $(B^+)^\S$ of fixed points of $B^+$ under $\S$ is generated by the elements $\D_J$, for $J$ running through the spherical orbits of $I$ under $\S$. 

We fix four parameters $(a,b,c,d)\in R$ such that $d^2 - ad - bc = 0$ and consider an LK-family $(f_i)_{i\in I}$ associated with $(a,b,c,d)$ that satisfy $f_i = f_{\s(i)}\s_V$ for every $(i,\s)\in I\times \S$, so that the associated twisted LK-representation $\psi^\S : (B^+)^\S \to \EL(V^\S)$ is defined. 

The aim of this section is to compute and study the map $\psi^\S_{\D_J}$ for a spherical orbit $J$ of $I$ under $\S$. For sake of brevity, we set $\psi^\S_J := \psi^\S_{\D_J}$ and $\p^\S_J := \p^\S_{\D_J}$ for any spherical orbit $J$.

\begin{nota}\label{notations ABCD}Any spherical orbit $J$ of $I$ under $\S$ is of one of the following types :
\begin{enumerate}
\item type A : $J = \{i\}$, 
\item type B : $J = \{i,j\}$ with $m_{i,j} = 2$, 
\item type C : $J = \{i,j,k\}$ with $m_{i,j} = m_{j,k} = m_{k,i} = 2$,  
\item type D : $J = \{i,j\}$ with $m_{i,j} = 3$.
\end{enumerate}
\end{nota}

Notice that types B and D only occur when $\S$ is of order two, and type C only occurs when $\S$ is of order three. We will make a constant use of these types A to D in the following.   

\begin{nota}Let $J$ be a spherical orbit of $I$ under $\S$. 
\begin{itemize}
\item We denote by $\Theta_J$ the set $\{\a_i\mid i\in J\}$ ; it is an orbit of $\Phi^+$ under $\S$.
\item In case D, $\{\a_i + \a_j\}$ is also an orbit of $\Phi^+$ under $\S$, that we denote by $\Theta'_J$. 
\end{itemize} 
\end{nota}

\begin{nota}\label{def des formes twisted}Let $J$ be a spherical orbit of $I$ under $\S$. Since we are assuming that $f_i = f_{\s(i)}\s_V$ for every $i\in I$, the linear forms $f_i$, for $i \in J$, coincide on $V^\S$. Depending on $J$, we define the following linear forms on $V^\S$ : 
\begin{enumerate}
\item for type A : $f_J = f_i|_{V^\S}$, where $J = \{i\}$,
\item for type B : $f_J = df_i|_{V^\S}$, for any $i\in J$,
\item for type C : $f_J = d^2f_i|_{V^\S}$, for any $i\in J$,
\item for type D : $f_J = f_j(bc\Id+a\p_i)|_{V^\S}$ and $f'_J = cf_j\p_i|_{V^\S}$, where $J = \{i,j\}$. 
\end{enumerate} 
\end{nota}

Notice that the two linear forms for case D really only depend on $J$ (not on the choice of $i$ and $j$ in $J$) since $f_i|_{V^\S} = f_j|_{V^\S}$ by assumption and since $f_j\p_i = f_i\p_j$ by definition of an LK-family. 

\begin{lem}\label{calcul intermediaire}If $i,\, j\in I$ with $i\neq j$, then $f_j(bc\Id+a\p_i) = f_j\p_i^2$ in $V^\star$. Moreover if $m_{i,j} = 3$, then $f_j\p_i^2 = f_i\p_j\p_i$ in $V^\star$.
\end{lem}
\proof By \cite[Lem. 11]{Ca}, the endomorphism $\p_i^2 - a\p_i -bc\Id$ has its image included in $Re_{\a_i}$, which is included in $\ker(f_j)$ if $i \neq j$ since then $f_j(e_{\a_i}) = 0$ by definition of an LK-family. So $f_j(\p_i^2 - a\p_i -bc\Id)$ is the zero form and this gives the first point. The second point is clear since then $f_j\p_i = f_i\p_j$ by definition of an LK-family. \qed 

\begin{pro}\label{forme de psiG}Let $J$ be a spherical orbit of $I$ under $\S$. Then  
\begin{enumerate}
\item $\psi^\S_{J} = \p^\S_{J}+f_J\bt e_{\Theta_J}$ if $J$ is of type A, B or C, 
\item $\psi^\S_{J} = \p^\S_{J} + f_J\bt e_{\Theta_J} + f'_J\bt e_{\Theta'_J}$ if $J$ is of type D.
\end{enumerate}
\end{pro}
\proof Recall that for every $i\in I$, $\psi_i$ is given by $\psi_i = \p_i + f_i\bt e_{\a_i}$. The result for type A is trivial, since then $\psi^\S_J = \psi_i|_{V^\S}$. For type B, we have $\psi^\S_J = \psi_i\psi_j|_{V^\S}$. But in that case, we also have $\p_i(e_{\a_j}) = de_{\a_j}$ (by definition of $\p_i$) and $f_{i}(e_{\a_j}) = 0$ (by definition of an LK-family), thus we get, thanks to the formulas of \cite[Rem. 1]{Ca}, 
$$\psi_i\psi_j = \p_i\p_j + f_i\p_j\bt e_{\a_i} + df_j\bt e_{\a_j},$$ 
whence the result since $f_i\p_j = df_i$ by definition of an LK-family. Similarly for type C, we have $\psi^\S_J = \psi_i\psi_j\psi_k|_{V^\S}$ and since $\p_i(e_{\a_j}) = de_{\a_j}$, $\p_i(e_{\a_k}) = \p_j(e_{\a_k}) = de_{\a_k}$, and $f_{i}(e_{\a_j}) = f_i(e_{\a_k}) = f_j(e_{\a_k}) = 0$, we get, thanks to \cite[Rem. 1]{Ca},  
$$\psi_i\psi_j\psi_k = \p_i\p_j\p_k + f_i\p_j\p_k\bt e_{\a_i} + df_j\p_k\bt e_{\a_j} + d^2f_k\bt e_{\a_k},$$ 
whence the result since $f_i\p_j\p_k = df_i\p_k = d^2f_i$ and $f_j\p_k = df_j$ by definition of an LK-family. Finally for type D, we have $\psi^\S_J = \psi_i\psi_j\psi_i|_{V^\S}$ and $\p_i(e_{\a_j}) = ae_{\a_j}+ce_{\a_i+\a_j}$, $\p_i\p_j(e_{\a_i}) = bce_{\a_j}$, $f_{i}(e_{\a_j}) = f_j(e_{\a_i}) = 0$, and $f_i\p_j(e_{\a_i}) = 0$ (since $f_i\p_j = f_j\p_i$ and $\p_i(e_{\a_i}) = 0$), hence we get, thanks to \cite[Rem. 1]{Ca}, 
$$\psi_i\psi_j\psi_i = \p_i\p_j\p_i + f_i\p_j\p_i\bt e_{\a_i} + \left[bcf_i+af_j\p_i\right]\bt e_{\a_j}+cf_j\p_i\bt e_{\a_i+\a_j},$$
whence the result by lemma \ref{calcul intermediaire} and by the fact that $f_i|_{V^\S} = f_j|_{V^\S}$. \qed

\begin{rem}\label{phiSigma est diag par blocks}Let $J$ be a spherical orbit of $I$ under $\S$ and let $\Theta$ be an orbit of $\Phi^+$ under $\S$. Recall that we denote by $M_{J}(\Theta)$ the $J$-mesh of $\Theta$, \ie the set $W_{J}(\Theta)\bigcap \Phi^+$. Notice that this subset of $\Phi^+$ is a union of orbits of $\Phi^+$ under $\S$.

It is clear that all the maps $\p_i$, for $i\in J$, stabilize the submodule of $V$ generated by the elements $e_\b$ for $\b$ running through $M_{J}(\Theta)$. So the map $\p^\S_{J}$ stabilizes the submodule of $V^\S$ generated by the elements $e_\Theta$ for $\Theta$ running through the orbits included in $M_{J}(\Theta)$, and hence the matrix of $\p^\S_{J}$ in the basis $(e_\Theta)_{\Theta \in \Phi^+/\S}$ of $V^\S$ is block diagonal, relatively to the partition of $\Phi^+/\S$ induced by those $J$-meshes $M_{J}(\Theta)$ for $\Theta \in \Phi^+/\S$.
\end{rem}

In subsections \ref{J singleton} to \ref{J avec mij=3} below, we explicit the different possible blocks of $\p^\S_{J}$, for the four possible types of $J$ and for the different possible $J$-meshes $M_{J}(\Theta)$ when $\Theta$ runs through $\Phi^+/\S$. For each such $J$-mesh $X = M_{J}(\Theta)$, we denote by $M_{X}$ the corresponding block in $\p^\S_J$ and by $P_{X}$ its characteristic polynomial.

\begin{nota}
For sake of brevity in what follows, we set $\check{d} := a- d$ in $R$. In particular, the roots of the polynomial $X^2 - a X - bc$ are $d$ and $\check{d}$, and $d\check{d} = -bc$.
\end{nota}

\subsection{Type A}\label{J singleton}\mbox{}\medskip 

We assume here that $J = \{i\}$, and hence that $\p^\S_J = \p_i|_{V^\S}$. 

\medskip

\noindent $\bullet$ Configuration {\bf A1} : $\Theta = \Theta_J$. The corresponding block is 
$$M_{\textrm{\bf A1}} = \begin{pmatrix}0\end{pmatrix}.$$

\noindent $\bullet$ Configuration {\bf A2} : 
\begin{picture}(35,15)(0,-4)
\put(20,0){\circle*{4}}
\put(25,0){\circle{10}}
\put(24,-3){\scriptsize $i$}
\put(0,-3){$\Theta$}
\put(15,-4){\dashbox(18,8){}}
\end{picture}  or 
\begin{picture}(55,0)(0,-4)
\put(20,0){\circle*{4}}
\put(25,0){\circle{10}}
\put(24,-3){\scriptsize $i$}
\put(40,0){\circle*{4}}
\put(45,0){\circle{10}}
\put(44,-3){\scriptsize $i$}
\put(0,-3){$\Theta$}
\put(15,-4){\dashbox(38,8){}}
\end{picture} or 
\begin{picture}(75,0)(0,-4)
\multiput(20,0)(20,0){3}{\circle*{4}}
\multiput(25,0)(20,0){3}{\circle{10}}
\multiput(24,-3)(20,0){3}{\scriptsize $i$}
\put(0,-3){$\Theta$}
\put(15,-4){\dashbox(58,8){}}
\end{picture}. 

The corresponding block is 
$$M_{\textrm{\bf A2}} = \begin{pmatrix}d\end{pmatrix}.$$

\noindent $\bullet$ Configuration {\bf A3} : 
\begin{picture}(35,25)(0,4)
\put(20,0){\circle*{4}}
\put(20,20){\circle*{4}}
\put(20,0){\line(0,1){20}}
\put(22,8){\scriptsize $i$}
\put(0,-3){\small  $\Theta_1$}
\put(0,17){\small $\Theta_2$}
\put(15,-4){\dashbox(10,8){}}
\put(15,16){\dashbox(10,8){}}
\end{picture}  or 
\begin{picture}(55,25)(0,4)
\put(20,0){\circle*{4}}
\put(20,20){\circle*{4}}
\put(20,0){\line(0,1){20}}
\put(22,8){\scriptsize $i$}
\put(0,-3){\small  $\Theta_1$}
\put(0,17){\small $\Theta_2$}
\put(15,-4){\dashbox(30,8){}}
\put(15,16){\dashbox(30,8){}}
\put(40,0){\circle*{4}}
\put(40,20){\circle*{4}}
\put(40,0){\line(0,1){20}}
\put(42,8){\scriptsize $i$}
\end{picture} or 
\begin{picture}(75,25)(0,4)
\multiput(20,0)(20,0){3}{\circle*{4}}
\multiput(20,20)(20,0){3}{\circle*{4}}
\multiput(20,0)(20,0){3}{\line(0,1){20}}
\multiput(22,8)(20,0){3}{\scriptsize $i$}
\put(0,-3){\small  $\Theta_1$}
\put(0,17){\small $\Theta_2$}
\put(15,-4){\dashbox(50,8){}}
\put(15,16){\dashbox(50,8){}}
\end{picture}. 

\bigskip

The corresponding block is 
$$M_{\textrm{\bf A3}} = \overset{\text{\scriptsize \ \ $e_{\Theta_1}$\  $e_{\Theta_2}$}}{\begin{pmatrix}a & b \\ c & 0 \end{pmatrix}},$$
whose characteristic polynomial is $$P_{\textrm{\bf A3}} = (X-d)(X-\check{d}).$$ 

\subsection{Type B}\label{J avec mij=2}\mbox{}\medskip 

We assume here that $J = \{i,j\}$ with $m_{i,j} = 2$, and hence that $\p^\S_{J} = (\p_i\p_j)|_{V^\S}$.


\medskip

\noindent $\bullet$ Configuration {\bf B1} : $\Theta = \Theta_J$. The corresponding block is 
$$M_{\textrm{\bf B1}} = \begin{pmatrix}0\end{pmatrix}.$$

\noindent $\bullet$ Configuration {\bf B2} : 
\begin{picture}(35,18)(0,-4)
\put(20,0){\circle*{4}}
\put(25,0){\circle{10}}
\put(24,-3){\scriptsize $i$}
\put(15,0){\circle{10}}
\put(12,-3){\scriptsize $j$}
\put(-3,-3){\small $\Theta$}
\put(8,-4){\dashbox(24,8){}}
\end{picture} or \ 
\begin{picture}(60,10)(0,-4)
\put(20,0){\circle*{4}}
\put(25,0){\circle{10}}
\put(24,-3){\scriptsize $i$}
\put(15,0){\circle{10}}
\put(12,-3){\scriptsize $j$}
\put(45,0){\circle*{4}}
\put(50,0){\circle{10}}
\put(49,-3){\scriptsize $i$}
\put(40,0){\circle{10}}
\put(37,-3){\scriptsize $j$}
\put(-3,-3){\small $\Theta$}
\put(8,-4){\dashbox(50,8){}}
\end{picture}. The corresponding block is 
$$M_{\textrm{\bf B2}} = \begin{pmatrix}d^2\end{pmatrix}.$$

\noindent $\bullet$ Configuration {\bf B3} : 
\begin{picture}(55,30)(-5,4)
\put(20,0){\circle*{4}}
\put(20,20){\circle*{4}}
\put(20,0){\line(0,1){20}}
\put(22,8){\scriptsize $i$}
\put(15,0){\circle{10}}
\put(12,-3){\scriptsize $j$}
\put(15,20){\circle{10}}
\put(12,17){\scriptsize $j$}
\put(-6,-3){\small  $\Theta_1$}
\put(-6,17){\small $\Theta_2$}
\put(8,-4){\dashbox(36,8){}}
\put(8,16){\dashbox(36,8){}}
\put(40,0){\circle*{4}}
\put(40,20){\circle*{4}}
\put(40,0){\line(0,1){20}}
\put(42,8){\scriptsize $j$}
\put(35,0){\circle{10}}
\put(32,-3){\scriptsize $i$}
\put(35,20){\circle{10}}
\put(32,17){\scriptsize $i$}
\end{picture}. The corresponding block is 

\smallskip

$$M_{\textrm{\bf B3}} = \overset{\text{\scriptsize \ \ $e_{\Theta_1}$\ \ $e_{\Theta_2}$}}{\begin{pmatrix}a d & b d \\ c d & 0 \end{pmatrix}},$$ 
whose characteristic polynomial is $$P_{\textrm{\bf B3}} = (X-d^2)(X-d\check{d}).$$ 

\noindent $\bullet$ Configuration {\bf B4} : 
\begin{picture}(80,40)(-30,14)
\put(20,0){\circle*{4}}
\put(0,20){\circle*{4}}
\put(40,20){\circle*{4}}
\put(20,40){\circle*{4}}
\put(20,0){\line(1,1){20}}
\put(20,0){\line(-1,1){20}}
\put(20,40){\line(1,-1){20}}
\put(20,40){\line(-1,-1){20}}
\put(33,5){\scriptsize $i$}
\put(33,30){\scriptsize $j$}
\put(2,5){\scriptsize $j$}
\put(5,30){\scriptsize $i$}
\put(-23,-3){\small  $\Theta_1$}
\put(-23,17){\small $\Theta_2$}
\put(-23,37){\small $\Theta_3$}
\put(15,-4){\dashbox(10,8){}}
\put(15,36){\dashbox(10,8){}}
\put(-5,16){\dashbox(50,8){}}
\end{picture}. The corresponding block is 

\medskip

$$M_{\textrm{\bf B4}} = \overset{\text{\scriptsize \ \ \ $e_{\Theta_1}$\ \ \ \ $e_{\Theta_2}$ \ \ $e_{\Theta_3}$}}{\begin{pmatrix}
a^2 & 2{ab} & b^2 \\
{ac} & {bc} & 0 \\
c^2 & 0 & 0
\end{pmatrix}},$$
whose characteristic polynomial is $$P_{\textrm{\bf B4}} = (X-d^2)(X-d\check{d})(X-\check{d}^2).$$ 

\noindent $\bullet$ Configuration {\bf B5} : 
\begin{picture}(120,40)(-30,14)
\put(20,0){\circle*{4}}
\put(0,20){\circle*{4}}
\put(40,20){\circle*{4}}
\put(20,40){\circle*{4}}
\put(20,0){\line(1,1){20}}
\put(20,0){\line(-1,1){20}}
\put(20,40){\line(1,-1){20}}
\put(20,40){\line(-1,-1){20}}
\put(33,8){\scriptsize $i$}
\put(33,29){\scriptsize $j$}
\put(3,7){\scriptsize $j$}
\put(5,30){\scriptsize $i$}
\put(-23,-3){\small  $\Theta_1$}
\put(-23,17){\small $\Theta_2$}
\put(-23,37){\small $\Theta_4$}
\put(15,-4){\dashbox(30,8){}}
\put(15,36){\dashbox(30,8){}}
\put(-5,16){\dashbox(30,8){}}
\put(40,0){\circle*{4}}
\put(20,20){\circle*{4}}
\put(60,20){\circle*{4}}
\put(40,40){\circle*{4}}
\put(40,0){\line(1,1){20}}
\put(40,0){\line(-1,1){20}}
\put(40,40){\line(1,-1){20}}
\put(40,40){\line(-1,-1){20}}
\put(53,7){\scriptsize $j$}
\put(52,30){\scriptsize $i$}
\put(35,16){\dashbox(30,8){}}
\put(73,17){\small $\Theta_3$}
\end{picture}. The corresponding block is 

\bigskip

$$M_{\textrm{\bf B5}} = \overset{\text{\scriptsize \ \ $e_{\Theta_1}$\ \ $e_{\Theta_2}$ \ \ $e_{\Theta_3}$  \  $e_{\Theta_4}$}}{\begin{pmatrix}
a^2 & {ab} & {ab} & b^2 \\
{ac} & 0 &{bc} & 0\\
{ac} & {bc} & 0 & 0 \\
c^2 & 0 & 0 & 0
\end{pmatrix}},$$
whose characteristic polynomial is $$P_{\textrm{\bf B5}} = (X-d^2)(X-d\check{d})^2(X-\check{d}^2).$$ 
Moreover, it is easily checked that the polynomial $P_{\textrm{\bf B4}}$ already annihilates $M_{\textrm{\bf B5}}$.

\subsection{Type C}\label{J  de card 3}\mbox{}\medskip

We assume here that $J = \{i,j,k\}$ with $m_{i,j} = m_{j,k} = m_{k,i} = 2$, and hence that $\p^\S_{J} = (\p_i\p_j\p_k)|_{V^\S}$.


\medskip

\noindent $\bullet$ Configuration {\bf C1} : $\Theta = \Theta_J$. The corresponding block is 
$$M_{\textrm{\bf C1}} = \begin{pmatrix}0\end{pmatrix}.$$

\noindent $\bullet$ Configuration {\bf C2} : \ 
\begin{picture}(30,15)(0,-2)
\put(20,0){\circle*{4}}
\put(25,-3){\circle{10}}
\put(24,-5){\scriptsize $i$}
\put(15,-3){\circle{10}}
\put(12,-5){\scriptsize $j$}
\put(20,6){\circle{10}}
\put(18,4){\scriptsize $k$}
\put(-3,-3){\small $\Theta$}
\put(8,-4){\dashbox(24,8){}}
\end{picture} \  or \ 
\begin{picture}(85,15)(0,-2)
\multiput(20,0)(25,0){3}{\circle*{4}}
\multiput(25,-3)(25,0){3}{\circle{10}}
\multiput(24,-5)(25,0){3}{\scriptsize $i$}
\multiput(15,-3)(25,0){3}{\circle{10}}
\multiput(12,-5)(25,0){3}{\scriptsize $j$}
\multiput(20,6)(25,0){3}{\circle{10}}
\multiput(18,4)(25,0){3}{\scriptsize $k$}
\put(-3,-3){\small $\Theta$}
\put(8,-4){\dashbox(75,8){}}
\end{picture}. The corresponding block is 

\smallskip

$$M_{\textrm{\bf C2}} = \begin{pmatrix}d^3\end{pmatrix}.$$

\noindent $\bullet$ Configuration {\bf C3} : \ 
\begin{picture}(95,30)(-5,4)
\multiput(20,0)(25,0){3}{\circle*{4}}
\multiput(20,20)(25,0){3}{\circle*{4}}
\multiput(20,0)(25,0){3}{\line(0,1){20}}
\multiput(15,0)(25,0){3}{\circle{10}}
\multiput(15,20)(25,0){3}{\circle{10}}
\multiput(25,0)(25,0){3}{\circle{10}}
\multiput(25,20)(25,0){3}{\circle{10}}
\put(-6,-3){\small  $\Theta_1$}
\put(-6,17){\small $\Theta_2$}
\put(8,-4){\dashbox(75,8){}}
\put(8,16){\dashbox(75,8){}}
\put(12,19){\scriptsize $j$}
\put(12,-3){\scriptsize $j$}
\put(23,19){\scriptsize $k$}
\put(23,-3){\scriptsize $k$}
\put(22,8){\scriptsize $i$}

\put(37,19){\scriptsize $i$}
\put(37,-3){\scriptsize $i$}
\put(48,19){\scriptsize $k$}
\put(48,-3){\scriptsize $k$}
\put(47,8){\scriptsize $j$}

\put(62,19){\scriptsize $j$}
\put(62,-3){\scriptsize $j$}
\put(73,19){\scriptsize $i$}
\put(73,-3){\scriptsize $i$}
\put(72,8){\scriptsize $k$}
\end{picture}. The corresponding block is 

\smallskip

$$M_{\textrm{\bf C3}} = \overset{\text{\scriptsize \ \ $e_{\Theta_1}$\ \ \ $e_{\Theta_2}$}}{\begin{pmatrix}a d^2 & b d^2 \\ c d^2 & 0 \end{pmatrix}},$$ 

whose characteristic polynomial is $$P_{\textrm{\bf C3}} = (X-d^3)(X-d^2\check{d}).$$ 

\noindent $\bullet$ Configuration {\bf C4} : \ 
\begin{picture}(140,40)(-35,-5)
\multiput(30,-20)(20,0){3}{\circle*{4}}
\multiput(0,0)(20,0){6}{\circle*{4}}
\multiput(30,20)(20,0){3}{\circle*{4}}
\multiput(30,-25)(20,0){3}{\circle{10}}
\multiput(30,25)(20,0){3}{\circle{10}}
\multiput(-5,0)(20,0){3}{\circle{10}}
\multiput(65,0)(20,0){3}{\circle{10}}
\multiput(30,-20)(20,0){3}{\line(-3,2){30}}
\multiput(60,0)(20,0){3}{\line(-3,2){30}}
\multiput(30,-20)(20,0){3}{\line(3,2){30}}
\multiput(0,0)(20,0){3}{\line(3,2){30}}
\put(28,-28){\scriptsize $k$}
\put(48,-28){\scriptsize $i$}
\put(68,-28){\scriptsize $j$}
\put(28,24){\scriptsize $k$}
\put(48,24){\scriptsize $i$}
\put(68,24){\scriptsize $j$}
\put(-8,-2){\scriptsize $k$}
\put(12,-2){\scriptsize $i$}
\put(32,-2){\scriptsize $j$}
\put(64,-2){\scriptsize $k$}
\put(84,-2){\scriptsize $i$}
\put(104,-2){\scriptsize $j$}

\put(48,8){\scriptsize $i$}
\put(66,10){\scriptsize $j$}
\put(86,10){\scriptsize $k$}

\put(10,10){\scriptsize $j$}
\put(30,10){\scriptsize $k$}

\put(39,-10){\scriptsize $j$}
\put(59,-11){\scriptsize $k$}
\put(86,-14){\scriptsize $i$}

\put(10,-14){\scriptsize $i$}
\put(-30,-24){$\Theta_1$}
\put(-30,17){$\Theta_4$}
\put(20,-23){\dashbox(60,8){}}
\put(20,16){\dashbox(60,8){}}
\put(-30,-3){$\Theta_2$}
\put(120,-3){$\Theta_3$}
\put(-15,-4){\dashbox(60,8){}}
\put(55,-4){\dashbox(60,8){}}
\end{picture}. 

\vspace{1cm}

The corresponding block is 
$$M_{\textrm{\bf C4}} = \overset{\text{\scriptsize \ \ $e_{\Theta_1}$\ \ \ $e_{\Theta_2}$ \ \ \ $e_{\Theta_3}$ \ \ \  $e_{\Theta_4}$}}{\begin{pmatrix}
a^2d & {abd} & {abd} & b^2d \\
{acd} & 0 &{bcd} & 0\\
{acd} & {bcd} & 0 & 0 \\
c^2d & 0 & 0 & 0
\end{pmatrix}},$$ 
whose characteristic polynomial is 
$$
P_{\textrm{\bf C4}} = (X-d^3)(X-d^2\check{d})^2(X-d\check{d}^2).
$$ 

Moreover, one can check that the polynomial $(X-d^3)(X-d^2\check{d})(X-d\check{d}^2)$ already annihilates $M_{\textrm{\bf C4}}$.

\medskip

\noindent $\bullet$ Configuration {\bf C5} : \ 
\begin{picture}(50,40)(-35,24)
\multiput(0,20)(0,20){2}{\circle*{4}}
\multiput(40,20)(0,20){2}{\circle*{4}}
\multiput(20,0)(0,20){4}{\circle*{4}}
\multiput(0,20)(40,0){2}{\line(0,1){20}}
\multiput(20,0)(0,40){2}{\line(0,1){20}}
\multiput(20,0)(20,20){2}{\line(-1,1){20}}
\multiput(20,0)(-20,20){2}{\line(1,1){20}}
\multiput(20,60)(-20,-20){2}{\line(1,-1){20}}
\multiput(20,60)(20,-20){2}{\line(-1,-1){20}}
\multiput(15,9)(0,39){2}{\scriptsize $j$}
\multiput(-5,28)(47,0){2}{\scriptsize $j$}
\put(5,6){\scriptsize $i$}
\put(32,50){\scriptsize $i$}
\put(5,50){\scriptsize $k$}
\put(32,6){\scriptsize $k$}
\put(10,22){\scriptsize $i$}
\put(26,34){\scriptsize $i$}
\put(26,22){\scriptsize $k$}
\put(10,34){\scriptsize $k$}
\put(-30,-3){$\Theta_1$}
\put(-30,57){$\Theta_4$}
\put(12,-4){\dashbox(16,8){}}
\put(12,56){\dashbox(16,8){}}
\put(-30,17){$\Theta_2$}
\put(-30,37){$\Theta_3$}
\put(-8,16){\dashbox(56,8){}}
\put(-8,36){\dashbox(56,8){}}
\end{picture}

\vspace{1.2cm}

The corresponding block is 
$$M_{\textrm{\bf C5}} = \overset{\text{\scriptsize \ \ $e_{\Theta_1}$\ \ \ $e_{\Theta_2}$ \ \ \ $e_{\Theta_3}$ \ \ \  $e_{\Theta_4}$}}{\begin{pmatrix}
a^3 & 3a^2b & 3ab^2 & b^3 \\
a^2c & 2abc & b^2c & 0\\
ac^2 & bc^2 & 0 & 0 \\
c^3 & 0 & 0 & 0
\end{pmatrix}},$$
whose characteristic polynomial is 
$$
P_{\textrm{\bf C5}} = (X-d^3)(X-d^2\check{d})(X-d\check{d}^2)(X-\check{d}^3).
$$ 

\medskip

\noindent $\bullet$ Configuration {\bf C6} : \ 
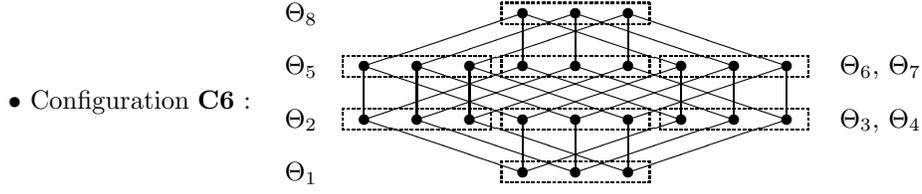
\begin{picture}(180,40)(-35,24)
\multiput(60,0)(20,0){3}{\circle*{4}}
\multiput(0,20)(20,0){9}{\circle*{4}}
\multiput(0,40)(20,0){9}{\circle*{4}}
\multiput(60,60)(20,0){3}{\circle*{4}}
\multiput(60,0)(20,0){3}{\line(-3,1){60}}
\multiput(60,20)(20,0){6}{\line(-3,1){60}}
\multiput(120,40)(20,0){3}{\line(-3,1){60}}
\multiput(60,0)(20,0){3}{\line(3,1){60}}
\multiput(0,20)(20,0){6}{\line(3,1){60}}
\multiput(0,40)(20,0){3}{\line(3,1){60}}
\multiput(60,0)(20,0){3}{\line(0,1){20}}
\multiput(60,40)(20,0){3}{\line(0,1){20}}
\multiput(0,20)(20,0){3}{\line(0,1){20}}
\multiput(120,20)(20,0){3}{\line(0,1){20}}
\put(-30,-3){$\Theta_1$}
\put(-30,57){$\Theta_8$}
\put(52,-4){\dashbox(56,8){}}
\put(52,56){\dashbox(56,8){}}
\put(-30,17){$\Theta_2$}
\put(180,17){$\Theta_3$, $\Theta_4$}
\multiput(-8,16)(60,0){3}{\dashbox(56,8){}}
\put(-30,37){$\Theta_5$}
\put(180,37){$\Theta_6$, $\Theta_7$}
\multiput(-8,36)(60,0){3}{\dashbox(56,8){}}
\end{picture}

\vspace{1.2cm}

\noindent(the labels $i$, $j$, $k$ on the edges are omitted).

The corresponding block is 
$$M_{\textrm{\bf C6}} = \overset{ \text{\scriptsize \ \ $e_{\Theta_1}$ \ \ $e_{\Theta_2}$ \ \ \ $e_{\Theta_3}$ \ \ \  $e_{\Theta_4}$ \ \ \ $e_{\Theta_5}$ \ \ \ $e_{\Theta_6}$ \ \ \ $e_{\Theta_7}$ \ \ \ $e_{\Theta_8}$ } }{\begin{pmatrix}
a^3 & a^2b & a^2b & a^2b & ab^2 & ab^2 & ab^2 & b^3 \\
a^2c & 0 & abc & abc & 0 & 0 & b^2c & 0\\
a^2c & abc & 0 & abc & 0 & b^2c & 0 & 0\\
a^2c & abc & abc & 0 & b^2c & 0 & 0 & 0\\
ac^2 & 0 & 0 & bc^2 & 0 & 0 & 0 & 0 \\
ac^2 & 0 & bc^2 & 0 & 0 & 0 & 0 & 0 \\
ac^2 & bc^2 & 0 & 0 & 0 & 0 & 0 & 0 \\
c^3 & 0 & 0 & 0 & 0 & 0 & 0 & 0
\end{pmatrix}},$$
whose characteristic polynomial is 
$$
P_{\textrm{\bf C6}} = (X-d^3)(X-d^2\check{d})^3(X-d\check{d}^2)^3(X-\check{d}^3).
$$ 
And one can check that the polynomial $P_{\textrm{\bf C5}}$ already annihilates $M_{\textrm{\bf C6}}$.

\subsection{Type D}\label{J avec mij=3}\mbox{}\medskip

We assume here that $J = \{i,j\}$ with $m_{i,j} = 3$, and hence that $\p^\S_{J} = (\p_i\p_j\p_i)|_{V^\S}$.


\medskip 

\noindent $\bullet$ Configuration {\bf D1} : $\Theta_J$ and $\Theta'_J$. The corresponding block is 
$$M_{\textrm{\bf D1}} = \begin{pmatrix}0 & 0\\ 0 & 0 \end{pmatrix}.$$

\noindent$\bullet$ Configuration {\bf D2} : \ 
\begin{picture}(30,15)(0,-2)
\put(20,0){\circle*{4}}
\put(25,0){\circle{10}}
\put(24,-3){\scriptsize $i$}
\put(15,0){\circle{10}}
\put(12,-3){\scriptsize $j$}
\put(-3,-3){\small $\Theta$}
\put(8,-4){\dashbox(24,8){}}
\end{picture} \ \  or\ \ \begin{picture}(60,15)(0,-2)
\put(20,0){\circle*{4}}
\put(25,0){\circle{10}}
\put(24,-3){\scriptsize $i$}
\put(15,0){\circle{10}}
\put(12,-3){\scriptsize $j$}
\put(45,0){\circle*{4}}
\put(50,0){\circle{10}}
\put(49,-3){\scriptsize $i$}
\put(40,0){\circle{10}}
\put(37,-3){\scriptsize $j$}
\put(-3,-3){\small $\Theta$}
\put(8,-4){\dashbox(50,8){}}
\end{picture}. The corresponding block is 
$$M_{\textrm{\bf D2}} = \begin{pmatrix}d^3\end{pmatrix}.$$


\noindent$\bullet$ Configuration {\bf D3} : \ \ 
\begin{picture}(50,35)(0,16)
\put(20,0){\circle*{4}}
\put(20,20){\circle*{4}}
\put(20,40){\circle*{4}}
\put(20,0){\line(0,1){40}}
\put(22,8){\scriptsize $i$}
\put(22,28){\scriptsize $j$}
\put(15,0){\circle{10}}
\put(12,-2){\scriptsize $j$}
\put(15,40){\circle{10}}
\put(12,37){\scriptsize $i$}
\put(-5,-3){\small  $\Theta_1$}
\put(-5,17){\small $\Theta_2$}
\put(-5,37){\small $\Theta_3$}
\put(8,-4){\dashbox(36,8){}}
\put(8,36){\dashbox(36,8){}}
\put(15,16){\dashbox(30,8){}}
\put(40,0){\circle*{4}}
\put(40,20){\circle*{4}}
\put(40,40){\circle*{4}}
\put(40,0){\line(0,1){40}}
\put(42,8){\scriptsize $j$}
\put(35,0){\circle{10}}
\put(32,-3){\scriptsize $i$}
\put(42,28){\scriptsize $i$}
\put(35,40){\circle{10}}
\put(32,38){\scriptsize $j$}
\end{picture}. The corresponding block is 

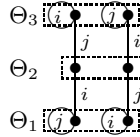
\vspace{4mm}

$$M_{\textrm{\bf D3}} = \overset{\text{\scriptsize $\overset{\mbox{}}{e_{\Theta_1}}$\ \ \ $e_{\Theta_2}$ \ \ \  $e_{\Theta_3}$}}{\begin{pmatrix}
{ad}^2 & {abd} & b^2 d \\
{acd} & {bcd} & 0 \\
c^2d & 0 & 0
\end{pmatrix}},$$
whose characteristic polynomial is 
$$
P_{\textrm{\bf D3}} = (X-d^3)(X^2+(d\check{d})^3).
$$


\noindent$\bullet$ Configuration {\bf D4} : \ \ \begin{picture}(70,45)(-20,22)
\put(20,0){\circle*{4}}
\put(0,20){\circle*{4}}
\put(40,20){\circle*{4}}
\put(0,40){\circle*{4}}
\put(40,40){\circle*{4}}
\put(20,60){\circle*{4}}
\put(20,0){\line(1,1){20}}
\put(20,0){\line(-1,1){20}}
\put(0,20){\line(0,1){20}}
\put(40,20){\line(0,1){20}}
\put(20,60){\line(1,-1){20}}
\put(20,60){\line(-1,-1){20}}
\put(33,5){\scriptsize $i$}
\put(33,50){\scriptsize $i$}
\put(43,28){\scriptsize $j$}
\put(-5,28){\scriptsize $i$}
\put(4,5){\scriptsize $j$}
\put(5,50){\scriptsize $j$}
\put(-23,-3){\small  $\Theta_1$}
\put(-23,17){\small $\Theta_2$}
\put(-23,37){\small $\Theta_3$}
\put(-23,57){\small $\Theta_4$}
\put(15,-4){\dashbox(10,8){}}
\put(15,56){\dashbox(10,8){}}
\put(-5,36){\dashbox(50,8){}}
\put(-5,16){\dashbox(50,8){}}
\end{picture}. The corresponding block is 

\vspace{4mm}

$$M_{\textrm{\bf D4}} = \overset{\text{\scriptsize \ \ \ \ \  \ \  $\overset{\mbox{}}{e_{\Theta_1}}$ \ \ \ \ \ \ \ \ \ $e_{\Theta_2}$  \ \ \ \ \   $e_{\Theta_3}$  \ \ \   $e_{\Theta_4}$}}{\begin{pmatrix}
a(a^2 + {bc}) & 2a^2b & 2{ab}^2 & b^3 \\
a^2c & 2{abc} & b^2c & 0\\
{ac}^2 & {bc}^2 & 0 & 0 \\
{c}^3 & 0 & 0 & 0
\end{pmatrix}},$$
whose characteristic polynomial is 
$$
P_{\textrm{\bf D4}} = (X-d^3)(X^2+(d\check{d})^3)(X-\check{d}^3).
$$


\noindent$\bullet$ Configuration {\bf D5} : \ \ \begin{picture}(110,55)(-20,20)
\put(20,0){\circle*{4}}
\put(0,20){\circle*{4}}
\put(40,20){\circle*{4}}
\put(0,40){\circle*{4}}
\put(40,40){\circle*{4}}
\put(20,60){\circle*{4}}
\put(20,0){\line(1,1){20}}
\put(20,0){\line(-1,1){20}}
\put(0,20){\line(0,1){20}}
\put(40,20){\line(0,1){20}}
\put(20,60){\line(1,-1){20}}
\put(20,60){\line(-1,-1){20}}
\put(33,8){\scriptsize $j$}
\put(33,49){\scriptsize $j$}
\put(33,27){\scriptsize $i$}
\put(23,27){\scriptsize $i$}
\put(3,7){\scriptsize $i$}
\put(3,50){\scriptsize $i$}
\put(-4,28){\scriptsize $j$}
\put(-23,-3){\small  $\Theta_1$}
\put(-23,17){\small $\Theta_2$}
\put(-23,37){\small $\Theta_4$}
\put(-23,57){\small $\Theta_6$}
\put(15,-4){\dashbox(30,8){}}
\put(15,56){\dashbox(30,8){}}
\put(-5,36){\dashbox(30,8){}}
\put(35,36){\dashbox(30,8){}}
\put(-5,16){\dashbox(30,8){}}
\put(40,0){\circle*{4}}
\put(20,20){\circle*{4}}
\put(60,20){\circle*{4}}
\put(20,40){\circle*{4}}
\put(60,40){\circle*{4}}
\put(40,60){\circle*{4}}
\put(40,0){\line(1,1){20}}
\put(40,0){\line(-1,1){20}}
\put(20,20){\line(0,1){20}}
\put(60,20){\line(0,1){20}}
\put(40,60){\line(1,-1){20}}
\put(40,60){\line(-1,-1){20}}
\put(53,7){\scriptsize $i$}
\put(53,50){\scriptsize $i$}
\put(62,28){\scriptsize $j$}
\put(35,16){\dashbox(30,8){}}
\put(73,17){\small $\Theta_3$}
\put(73,37){\small $\Theta_5$}
\end{picture}. The corresponding block is 

\vspace{5mm}

$$M_{\textrm{\bf D5}} = \overset{\text{\scriptsize \ \ \ \ \ \ \ $\overset{\mbox{}}{e_{\Theta_1}}$ \ \ \ \ \ \ \ \ $e_{\Theta_2}$ \ \ \ \ $e_{\Theta_3}$ \ \ \ $e_{\Theta_4}$  \ \ \  $e_{\Theta_5}$ \ \  $e_{\Theta_6}$}}{\begin{pmatrix}
a(a^2 + {bc}) & a^2b & a^2b & {ab}^2 & {ab}^2 & b^3 \\
a^2c & {abc} & {abc} & 0 & b^2c & 0\\
a^2c & {abc} & {abc} & b^2c & 0 & 0\\
{ac}^2 & 0 & {bc}^2 & 0 & 0 & 0 \\
{ac}^2 & {bc}^2 & 0 & 0 & 0 & 0\\
{c}^3 & 0 & 0 & 0 & 0 & 0
\end{pmatrix}},$$
whose characteristic polynomial is 
$$
P_{\textrm{\bf D5}} = (X-d^3)(X^2+(d\check{d})^3)^2(X-\check{d}^3).
$$
Moreover, it is easily checked that the polynomial $P_{\textrm{\bf D4}}$ already annihilates $M_{\textrm{\bf D5}}$.


\subsection{Consequences on annihilating polynomials}\label{Properties twisted}\mbox{}\medskip

\begin{lem}\label{puissance des psii}
Let $\p \in \EL(V)$, $f,\, f' \in V^\star$ and $e,\, e'\in V$ be such that $\p(e) = \p(e') = 0$ and $f(e') = f'(e) = 0$. We set $f_e = f(e)$ and $f'_{e'} = f'(e')$. Then we get the following identities in $\EL(V)$, for all $n \in \N$ :
\begin{enumerate}
\item $\displaystyle (\p + f\bt e)^n = \p^n + f\left[\sum_{k = 0}^{n-1}f_e^k\p^{n-1-k}\right]\bt e$,
\item $\displaystyle (\p + f\bt e + f'\bt e')^n = \p^n + f\left[\sum_{k = 0}^{n-1}f_e^k\p^{n-1-k}\right]\bt e + f'\left[\sum_{k = 0}^{n-1}f_{e'}'^k\p^{n-1-k}\right]\bt e'$.
\end{enumerate}
\end{lem}
\proof By induction on $n$, using the identities of \cite[Rem. 1]{Ca}. \qed

\begin{lem}\label{calcul intermediaire 2}Let $J = \{i,j\}$ be an orbit of $I$ under $\S$, with $m_{i,j} = 3$. Then $f_J(e_{\Theta'_J}) = f'_J(e_{\Theta_J}) = 0$ and $f_J(e_{\Theta_J}) = f'_J(e_{\Theta'_J})$.
\end{lem}
\proof With the formulas of notation \ref{def des formes twisted} and lemma \ref{calcul intermediaire}, we get for the first point :

$f_J(e_{\Theta'_J}) = f_i\p_j\p_i(e_{\a_i+\a_j}) = f_i\p_j(be_{\a_j}) = 0$, and 

$f'_J(e_{\Theta_J}) = cf_i\p_j(e_{\a_i}+e_{\a_j}) = cf_j\p_i(e_{\a_i}) + cf_i\p_j(e_{\a_j}) = 0$, 

\noindent and for the second point :

$f_J(e_{\Theta_J}) = f_i\p_j\p_i(e_{\a_i}+e_{\a_j}) = f_i\p_j\p_i(e_{\a_j}) = f_i(bce_{\a_i}) = bcf_i(e_{\a_i})$, and 

$f'_J(e_{\Theta'_J}) = cf_i\p_j(e_{\a_i+\a_j}) = cf_i(be_{\a_i}) = bcf_i(e_{\a_i})$. \qed

\begin{pro}\label{forme des P de PsiSigma}Let $J$ be a spherical orbit of $I$ under $\S$, and consider $P \in R[X]$. Then there exists a polynomial $Q \in R[X]$ (that depends on $J$ and $P$) such that : 
\begin{enumerate}
\item $P(\psi^\S_J) = P(\p^\S_J) + f_JQ(\p^\S_J)\bt e_{\Theta_J}$ if $J$ is of type A, B or C,
\item $P(\psi^\S_J) = P(\p^\S_J) + f_JQ(\p^\S_J)\bt e_{\Theta_J} + f'_JQ(\p^\S_J)\bt e_{\Theta'_J}$ if $J$ is of type D.
\end{enumerate}
\end{pro}
\proof For types A, B and C, the result follows from proposition \ref{forme de psiG} and lemma \ref{puissance des psii} (i), which applies since $\p^\S_{J}(e_{\Theta_J}) = 0$ (\cf configurations {\bf A1}, {\bf B1} and {\bf C1}). For type D, lemma \ref{puissance des psii} (ii) applies to the decomposition of $\psi^\S_J$ given in proposition \ref{forme de psiG} (ii) since $\p^\S_{J}(e_{\Theta_J}) = \p^\S_{J}(e_{\Theta'_J}) = 0$ (\cf configuration {\bf D1}), and since $f_J(e_{\Theta'_J}) = f'_J(e_{\Theta_J}) = 0$ by the previous lemma. Hence we get two polynomials $Q$ and $Q'$ such that $P(\psi^\S_J) = P(\p^\S_J) + f_JQ(\p^\S_J)\bt e_{\Theta_J} + f'_JQ'(\p^\S_J)\bt e_{\Theta'_J}$. But the two polynomials $Q$ and $Q'$ thus obtained are in fact equal since, again by the previous lemma, $f_J(e_{\Theta_J}) = f'_J(e_{\Theta'_J})$. \qed 

\begin{rem}\label{coeffs de QJ}More precisely, if $P = \sum_{n=0}^dp_nX^n$ in proposition \ref{forme des P de PsiSigma}, then $Q = \sum_{n=0}^{d-1}q_nX^n$ with $q_{d-1} = p_d$ and $q_{n-1} = p_n+q_nf_J(e_{\Theta_J})$ for every $1 \leqslant n\leqslant d-1$.
\end{rem}

\begin{nota}\label{les pols PJ}Let $J$ be a spherical orbit of $I$ under $\S$. Depending on $J$, we define a polynomial $P_J \in R[X]$ by : 
\begin{enumerate}
\item $P_J = (X-d)(X-\check{d})$ if $J$ is of type A, 
\item $P_J = (X-d^2)(X-d\check{d})(X-\check{d}^2)$ if $J$ is of type B, 
\item $P_J = (X-d^3)(X-d^2\check{d})(X-d\check{d}^2)(X-\check{d}^3)$ if $J$ is of type C, 
\item $P_J = (X-d^3)(X^2+(d\check{d})^3)(X-\check{d}^3)$ if $J$ is of type D.
\end{enumerate}
\end{nota}


\begin{lem}\label{image des P(phiJ)}Let $J$ be a spherical orbit of $I$ under $\S$. Then 
\begin{enumerate}
\item $P_J(\p^\S_{J})(e_\Theta) = 0$ when $\Theta \neq \Theta_J$ (resp. when $\Theta \neq \Theta_J$ and $\Theta \neq \Theta'_J$) if $J$ is of type A, B or C (resp. D),  
\item $P_J(\p^\S_{J})(e_{\Theta}) = P_J(0)e_{\Theta}$ when $\Theta = \Theta_J$ (resp. when $\Theta = \Theta_J$ or $\Theta = \Theta'_J$) if $J$ is of type A, B or C (resp. D).
\end{enumerate}
\end{lem}
\proof In view of the results of subsections \ref{J singleton} to \ref{J avec mij=3} above, the polynomial $P_J$ annihilates the non-zero blocks of $\p^\S_{J}$, hence we get that $P_J(\p^\S_{J})(e_\Theta) = 0$ if $\Theta \neq \Theta_J$ (resp. if $\Theta \neq \Theta_J$ and $\Theta \neq \Theta'_J$) for types A, B and C (resp. D). Moreover since $\p^\S_{J}(e_{\Theta_J}) = 0$ for all cases, and since $\p^\S_{J}(e_{\Theta'_J}) = 0$ for type D, we get that $P_J(\p^\S_{J})(e_{\Theta_J}) = P_J(0)e_{\Theta_J}$ for all cases, and that $P_J(\p^\S_{J})(e_{\Theta'_J}) = P_J(0)e_{\Theta'_J}$ for type~D. \qed

\begin{pro}\label{image des P(psiJ)}Let $J$ be a spherical orbit of $I$ under $\S$. The image of $P_J(\psi^\S_J)$ is included in $Re_{\Theta_J}$ (resp. in $Re_{\Theta_J}\oplus Re_{\Theta'_J}$) if $J$ is of type A, B or C (resp. D). As a consequence, the polynomial $(X-f_J(e_{\Theta_J}))P_J$ annihilates $\psi^\S_J$.
\end{pro}
\proof The previous lemma shows that the image of $P_J(\p^\S_J)$ is included in $Re_{\Theta_J}$ (resp. in $Re_{\Theta_J}\oplus Re_{\Theta'_J}$) for types A, B and C (resp. D). But by proposition \ref{forme des P de PsiSigma}, $P_J(\psi^\S_J)$ is the sum of $P_J(\p^\S_J)$ and of some endomorphism of $V^\S$ whose image is included in $Re_{\Theta_J}$ (resp. in $Re_{\Theta_J}\oplus Re_{\Theta'_J}$) for types A, B and C (resp. D), whence the first point of the proposition. 

The second point for types A, B and C (resp. D) follows from proposition \ref{forme de psiG} (i) (resp. proposition \ref{forme de psiG} (ii) and lemma \ref{calcul intermediaire 2}), and the studies of configurations {\bf A1}, {\bf B1}, and {\bf C1} (resp. {\bf D1}), which show that $e_{\Theta_J}$ (resp. each of $e_{\Theta_J}$ and $e_{\Theta'_J}$) is an eigenvector of $\psi^\S_J$ for the eigenvalue $f_J(e_{\Theta_J})$. \qed

\begin{rem}Let us denote by $f$ the common value of the $f_i(e_{\a_i})$ for $i\in J$ (recall that we are assuming that $f_i = f_{\s(i)}\s_V$ for all $i\in I$). Then it is easily checked on the definitions that the value of $f_J(e_{\Theta_J})$ is $f$ (resp. $df$, resp. $d^2f$, resp. $bcf$) if $J$ is of type A (resp. B, resp. C, resp. D).
\end{rem}


\section{The spherical case}\label{The spherical case}

We assume here that $\G$ is of spherical type $A_{n}$ ($n\geqslant2$), $D_{n}$ ($n\geqslant 4$) or $E_6$. 

In all the cases but $D_4$, we fix $\S = \Aut(\G)$ (of order two), and when $\G$ is of type $D_4$, we take for $\S$ a subgroup of $\Aut(\G)$ (which is of order six) of order two or three. The Coxeter matrix $\G^\S = (m^\S_{J,K})_{J,K\in I^\S}$ --- which encodes the presentations of the Coxeter group $W^\S$, the Artin-Tits monoid $(B^+)^\S$ and the Artin-Tits group $B^\S$ --- is then of type $B_{n}$ if $\G = A_{2n-1}$, $A_{2n}$ or $D_{n+1}$ with $|\S| = 2$, of type $F_4$ if $\G = E_6$, or of type $G_2$ if $\G = D_4$ with $|\S| = 3$ (see for example \cite{C2,Ca1}).

\medskip

We fix an integral domain $R$ and $(a,b,c,d)\in R^4$ such that $d^2 - ad - bc = 0$. Let $(f_i)_{i\in I} \in (V^\star)^I$ be an LK-family relatively to $(a,b,c,d)$ and consider the associated LK-representation $\psi : B^+\to \EL(V)$ of $B^+$. 

We will always assume in this section that the elements $b$, $c$ and $d$ are non-zero in $R$, so that they become invertible in the field of fractions $\K$ of $R$. By extension of the scalars from $R$ to $\K$, we will consider $\psi$ as acting on the $\K$-vector space $V_\K = V\otimes_R\K$. It then follows from proposition \ref{cas favorables} that $\psi$ induces a twisted LK-representation $\psi^\S : (B^+)^\S \to \EL(V_\K^\S)$ of $(B^+)^\S$ over $V^\S_\K = \K\otimes_RV^\S$.


By \cite[Section 3.2]{Ca}, the LK-family $(f_i)_{i\in I}$, seen as an element of $(V_\K^\star)^I$, is entirely determined by the common value $f$ of the elements $f_i(e_{\a_i})$, $i\in I$. Moreover when $f \neq 0$, then $\psi$ has its image included in $\GL(V_\K)$ and hence induces an LK-representation $\psi_{\gr} : B\to \GL(V_\K)$ of the Artin-Tits group $B$, and a twisted LK-representation $\psi_{\gr}^\S : B^\S \to \GL(V_\K^\S)$ of the Artin-Tits group $B^\S$. 

\medskip

In view of the results of the previous section, in order to really understand the maps $\psi^\S_J = \psi^\S_{\D_J}$, $J\in I^\S$, we need to list the possible configurations of the $J$-meshes $M_{J}(\Theta)$ in $\Phi^+$, when $\Theta$ runs through $\Phi^+/\S$. This is done in the following lemma.

\begin{lem}\label{nombre de confs}Fix an orbit $J$ of $I$ under $\S$. We list, in the following tables, the number $N_{\textrm{\bf X}}$ of occurrences of configuration {\bf X} (with the notations of subsections \ref{J singleton} to \ref{J avec mij=3}) among the $J$-meshes $M_{J}(\Theta)$ in $\Phi^+$, when $\Theta$ runs through $\Phi^+/\S$. (The configurations that do not occur are omitted).
%
%
\begin{itemize}
\item If $\G = A_{2n-1}$, or $D_{n}$ with $|\S|=2$, or $E_6$, then $J$ can be of type $A$ or $B$ and we get :
\end{itemize}
\begin{center}
\begin{tabular}{c|c|c|c||c|c|c|c|}
\cline{2-8}
 & \multicolumn{3}{c||}{if $J$ is of type A} & \multicolumn{4}{c|}{if $J$ is of type B}\\ 
 & $N_{\textrm{\bf A1}}$ & $N_{\textrm{\bf A2}}$ & $N_{\textrm{\bf A3}}$ & $N_{\textrm{\bf B1}}$ & $N_{\textrm{\bf B2}}$ & $N_{\textrm{\bf B3}}$ &  $N_{\textrm{\bf B4}}$  \\ \hline
\multicolumn{1}{|c||}{$A_{2n-1}$} & $1$ & $(n-1)^2$ & $n-1$ & $1$ & $(n-2)^2$ & $2n-4$ & $1$ \\ \hline
\multicolumn{1}{|c||}{$D_{n}$ ($|\S|\! =\! 2$)} & $1$ & $n^2\!-\!6n\!+\!10$ & $2n\!-\!5$ & $1$ & $(n\!-\!2)(n\!-\!3)$ & $0$ & $n\!-\!2$ \\ \hline
\multicolumn{1}{|c||}{$E_6$} & $1$ & $9$ & $7$ & $1$ & $6$ & $4$ & $3$ \\ \hline
\end{tabular}
\end{center}

\begin{itemize}
\item If $\G = A_{2n}$, then $J$ can be of type $B$ (if $n \geqslant 2$) or $D$ and we get :
\end{itemize}
\begin{center}
\begin{tabular}{|c|c|c|c||c|c|c|} 
\hline
 \multicolumn{4}{|c||}{if $J$ is of type B} & \multicolumn{3}{c|}{if $J$ is of type D}\\ 
 $N_{\textrm{\bf B1}}$ & $N_{\textrm{\bf B2}}$ & $N_{\textrm{\bf B3}}$ & $N_{\textrm{\bf B4}}$ & $N_{\textrm{\bf D1}}$ & $N_{\textrm{\bf D2}}$ &  $N_{\textrm{\bf D3}}$  \\ \hline
 $1$ & $(n-1)(n-2)$ & $2n-3$ & $1$ & $1$ & $(n-1)^2$ & $n-1$ \\ \hline
\end{tabular}
\end{center}

\begin{itemize}
\item If $\G = D_{4}$ with $|\S| = 3$, then $J$ can be of type $A$ or $C$ and we get :
\end{itemize}

\begin{center}
\begin{tabular}{|c|c|c||c|c|c|}
\hline
\multicolumn{3}{|c||}{if $J$ is of type A} & \multicolumn{3}{c|}{if $J$ is of type C}\\ 
$N_{\textrm{\bf A1}}$ & $N_{\textrm{\bf A2}}$ & $N_{\textrm{\bf A3}}$ & $N_{\textrm{\bf C1}}$ & $N_{\textrm{\bf C2}}$ & $N_{\textrm{\bf C5}}$  \\ \hline
 $1$ & $1$ & $2$ & $1$ & $1$ & $1$ \\ \hline
\end{tabular}
\end{center}
\end{lem}
\proof Left to the reader. \qed

\subsection{Irreducibility}\label{Irreducibility}\mbox{}\medskip

In this subsection, we exclude the case $\G = A_{2n}$, so that there is no orbit of type D in $I$, and hence no configurations {\bf D1} to {\bf D5} in $\Phi^+$.

\begin{lem}\label{on peut descendre}
For all $\Theta \in \Phi^+/\S$ with $\dep(\Theta)\geqslant 2$, there exists $J\in I^\S$  such that $\dep(r_J(\Theta)) < \dep(\Theta)$.
\end{lem}
\proof Fix $\a \in \Theta$ and $i\in I$ such that $\dep(s_i(\a)) <\dep(\a)$, and let us denote by $J$ the orbit of $i$ under $\S$. Necessarily, the $J$-mesh $M_J(\Theta)$ is in configuration {\bf A3}, {\bf B3}, {\bf B4} or {\bf C5}, and $\Theta$ is not of minimal depth in $M_J(\Theta)$.

In case of a configuration {\bf A3}, {\bf B3}, {\bf B4} with $\Theta$ of maximal depth in $M_J(\Theta)$, or {\bf C5} with $\Theta$ of maximal or second-maximal depth in $M_J(\Theta)$, then $\dep(r_J(\Theta)) < \dep(\Theta)$ and we are done. In the other cases, then ($\dep(r_J(\Theta)) \geqslant \dep(\Theta)$ and) a look at the graph of $\Phi^+$ in the different possible situations reveals that there exists $i' \in I$ with $\dep(s_{i'}(\a)) < \dep(\a)$ for which, if we denote by $J'$ its orbit under $\S$, then $M_{J'}(\Theta)$ is in configuration {\bf A3} or {\bf B3}, so here again we are done, up to changing $J$ for $J'$. \qed

\medskip

In the following proposition, we prove the first part of our irreducibility criterion.

\begin{nota}\label{def des QJ}
If $J$ if an orbit of $I$ under $\S$, we consider the polynomial $P_J$ as in notation \ref{les pols PJ}. In view of proposition \ref{forme des P de PsiSigma}, there exist a polynomial $Q_J\in R[X]$ such that $P_J(\psi^\S_J) = P_J(\p^\S_J) + f_JQ_J(\p^\S_J)\bt e_{\Theta_J}$. 

By lemma \ref{image des P(phiJ)}, we get $P_J(\psi^\S_J)(e_\Theta) = f_JQ_J(\p^\S_J)(e_\Theta) e_{\Theta_J}$ for every $\Theta\in \Phi^+/\S$ with $\Theta\neq\Theta_J$. Moreover it can be checked that $f_JQ_J(\p^\S_J)(e_{\Theta_K}) = 0$ for $J,\, K\in I^\S$ with $m_{J,K}^\S = 2$. 
\end{nota}

\begin{pro}\label{first part}
Assume that $f_JQ_J(\p^\S_J)(e_{\Theta_K}) \neq 0$ for every $J$, $K\in I^\S$ with $m^\S_{J,K}\geqslant 3$. If $U$ is an invariant subspace of $V_\K^\S$ under $\psi^\S$ such that $P_{K}(\psi^\S_{K})(U)\neq \{0\}$ for some $K \in I^\S$, then $U = V_\K^\S$.
\end{pro}
\proof By proposition \ref{image des P(psiJ)} and since $U$ is invariant, we get $P_J(\psi^\S_J)(U) \subseteq U \bigcap Re_{\Theta_J}$ for all $J \in I^\S$, whence $e_{\Theta_{K}} \in U$ since $P_{K}(\psi^\S_{K})(U)\neq \{0\}$. Now if $J \in I^\S$ is distinct from $K$, we get that $P_J(\psi^\S_J)(e_{\Theta_K}) = f_JQ_J(\psi^\S_J)(e_{\Theta_K})e_{\Theta_J}$ belongs to $U$, whence $e_{\Theta_J}\in U$ if $m^\S_{J,K} \geqslant 3$ since then, by assumption, $f_JQ_J(\p^\S_J)(e_{\Theta_K})\neq 0$. We thus get that $e_{\Theta_J} \in U$ for all $J \in I^\S$ by connectivity of $\G^\S$. 

Now consider $\Theta \in \Phi^+/\S$, of depth $p\geqslant 2$, such that every $e_{\Theta'}$ belongs to $U$ whenever $\dep(\Theta') < p$. By lemma \ref{on peut descendre}, there exists $J\in I^\S$ such that $\dep(r_J(\Theta)) < \dep(\Theta)$, so in particular $e_{r_J(\Theta)} \in U$. Necessarily, the $J$-mesh $M_J(\Theta)$ is in configuration {\bf A3}, {\bf B3}, {\bf B4} with $\Theta$ of maximal depth in $M_J(\Theta)$, or {\bf C5} with $\Theta$ of maximal or second-maximal depth in $M_J(\Theta)$, and hence $\psi^\S_J(e_{r_J(\Theta)})$ is a linear combination of $e_\Theta$, $e_{r_J(\Theta)}$, $e_{\Theta_J}$ and possibly some other $e_{\Theta'}$ for $\Theta' \subseteq M_J(\Theta)$ with $\dep(\Theta')<p$, in which the coefficient of $e_\Theta$ is $c$, $cd$, $c^2$, $c^3$ or $bc^2$ respectively. 

Since $U$ is invariant by $\psi^\S_J$ and since $bcd \neq 0$ by assumption, we thus get, in all the cases, the element $e_\Theta$ as a linear combination of elements of $U$, hence $e_\Theta \in U$ and we conclude by induction that $U = V_\K^\S$. \qed

\medskip

The following lemma explicits the condition on $f_JQ_J(\p^\S_J)(e_{\Theta_K})$ of the previous proposition (recall that $f$ is the common value of the $f_i(e_{\a_i})$ for $i\in I$, and that we set $\check{d} = a-d$) :

\begin{lem}\label{tech result 1}
Consider two orbits $J$, $K \in I^\S$ with $m^\S_{J,K}\geqslant 3$. Then the coefficient $f_JQ_J(\p^\S_J)(e_{\Theta_K})$ is equal to : 
\begin{enumerate}
\item $-|K|af$ if $J$ is of type A,
\item $ad^2f(\check{d}^2-df)$ if $J$ is of type B,
\item $ad^5f(-d^3f^2 + ad\check{d}^2f - \check{d}^5)$ if $J$ is of type C.
\end{enumerate}
\end{lem}
\proof First notice that the assumption $m^\S_{J,K}\geqslant 3$ implies that $m_{j,k} = 3$ for some $(j,k)\in J\times K$ (see \cite[Section 3]{C2} or \cite[Section 4.2]{Ca1}). Now it is easy to see that if $J$ is of type A (resp. B, resp. C), then $\Theta_K$ necessarily appears at the bottom of a configuration {\bf A3} (resp. {\bf B3} or {\bf B4}, resp. {\bf C5}). The result then follows from direct computations, using the expression of $Q_J$ deduced from remark \ref{coeffs de QJ} and notation \ref{les pols PJ}, and the expressions of the encountered $f_J(e_\Theta)$ as multiples of $f$ deduced from \cite[Table 1]{Ca}. \qed

\medskip

Now we turn to the second part of the irreducibility criterion.

\begin{pro}\label{second part}Assume that there exists a prime ideal $Q$ of $R$ with $bcd \not\in Q$ and $\Image(f_i) \subseteq Q$ for all $i\in I$. If $U$ is an invariant subspace of $V_\K^\S$ under $\psi^\S$ such that $P_{J}(\psi^\S_{J})(U) = \{0\}$ for all $J \in I^\S$, then $U = \{0\}$.
\end{pro}
\proof Recall that the image of $P_{J}(\psi^\S_{J})$ is included in $\K{e_{\Theta_J}}$. Since $U$ is invariant, we get $P_{J}(\psi^\S_{J})\psi^\S_g(U) = \{0\}$ for every $(J,g) \in I^\S\times(B^+)^\S$, and hence, if we denote by $L_{(J,g),\Theta}$ the element of $R$ such that $P_{J}(\psi^\S_{J})\psi^\S_g(e_\Theta) = L_{(J,g),\Theta}e_{\Theta_J}$, and by $L$ the matrix $(L_{(J,g),\Theta})_{(J,g)\in I^\S\times (B^+)^\S, \Theta \in \Phi^+/\S}$, we get that $U$ is included in $\ker(L)$. 

Now we claim that under the assumptions of the proposition, the matrix $L$ is non-singular, which thus implies that $U = \{0\}$. To prove our claim, it suffices to construct some pairs $(J_\Theta, g_\Theta) \in I^\S\times (B^+)^\S$ for $\Theta$ running through $\Phi^+/\S$, such that the square submatrix $L' = (L_{(J_\Theta,g_\Theta),\Omega})_{\Theta,\Omega \in \Phi^+/\S}$ of $L$ is invertible. We construct the pairs $(J_\Theta, g_\Theta) \in I^\S\times (B^+)^\S$ by induction on $\dep(\Theta)$ as follows. 

If $\dep(\Theta) = 1$, \ie if $\Theta = \Theta_J$ for some $J \in I^\S$, we set $J_\Theta = J$ and $g_\Theta = 1$. Now by induction if $\dep(\Theta)\geqslant 2$, then we fix some $J\in I^\S$ such that $\dep(r_J(\Theta))<\dep(\Theta)$ (which exists by lemma \ref{on peut descendre}) and we set $J_\Theta = J_{r_J(\Theta)}$ and $g_\Theta = g_{r_J(\Theta)}\D_J$. 

We are going to see that 
\begin{enumerate}
\item $P_{J_\Theta}(\p_{J_\Theta}^\S)\p^\S_{g_\Theta}(e_\Theta) = \pm b^kc^ld^me_{\Theta_{J_\Theta}}$ for some $k,\, l,\, m \in \N$,
\item $P_{J_\Theta}(\p_{J_\Theta}^\S)\p^\S_{g_\Theta}(e_\Omega) = 0$ for $\Omega \in \Phi^+/\S$ with $\dep(\Omega) \geqslant \dep(\Theta)$ and $\Omega \neq \Theta$.
\end{enumerate}

Since $P_{J}(\psi^\S_{J})\psi^\S_g(e_\Theta) \equiv P_{J}(\p^\S_{J})\p^\S_g(e_\Theta)\mod Q$, points (i) and (ii) will show that the matrix $L'$ is lower triangular modulo $Q$ --- if we choose an order on $\Phi^+/\S$ that is consistent with depth --- with non-zero diagonal coefficients modulo $Q$, since $bcd \not \in Q$. Hence this will show that the determinant of $L'$ is non-zero modulo $Q$ (recall that $Q$ is prime), and therefore is non-zero in $R$, whence the result. 

So it remains to prove (i) and (ii). 

Let us prove (i). If $\Theta = \Theta_J$ for some $J \in I^\S$, then by lemma \ref{image des P(phiJ)} we get $P_{J}(\p_{J}^\S)(e_{\Theta_J}) = P_J(0)e_{\Theta_J}$, with $P_J(0) = -bc$ (resp. $(bc)^3$, resp. $(bc)^6$) if $J$ is of type A (resp. B, resp. C), whence the result when $\dep(\Theta) = 1$. Now assume that $\dep(\Theta)\geqslant 2$ and consider the fixed $J\in I^\S$ such that $\dep(r_J(\Theta))<\dep(\Theta)$. 

If $M_J(\Theta)$ is in configuration {\bf A3}, {\bf B3}, {\bf B4} or {\bf C5}, with $\Theta$ of maximal depth in $M_J(\Theta)$, then $\p_J^\S(e_\Theta) = \lambda e_{r_J(\Theta)}$ with $\lambda = b$, $bd$, $b^2$ or $b^3$ respectively, whence the result by induction. The only other possible configuration is {\bf C5} with $\Theta$ of second maximal depth in $M_J(\Theta)$, which occurs when $\G = D_4$ with $|\S| = 3$. If we label the vertices of $D_4$ as in \cite[Planche IV]{Bo}, we thus get $I^\S = \{J,K\}$ with $J = \{1,3,4\}$ and $K = \{2\}$, $\Theta = \{\a_1+\a_2+\a_3, \a_1+\a_2+\a_4, \a_2+\a_3+\a_4\}$, $g_\Theta = \D_K\D_J$, and $J_\Theta = J$, whence $\p^\S_{g_\Theta}(e_\Theta) = b^3ce_{\Theta_J}$ and $P_{J_\Theta}(\p_{J_\Theta}^\S)\p^\S_{g_\Theta}(e_{\Theta}) = b^9c^7e_{\Theta_{J_\Theta}}$. 

Let us prove (ii). By lemma \ref{image des P(phiJ)}, it suffices to see that $\p^\S_{g_\Theta}(e_\Omega)$ is a linear combination of some $e_{\Theta'}$ with $\Theta'\neq \Theta_{J_\Theta}$. This is clear if $\Theta = \Theta_J$ since then $g_\Theta = 1$. So assume that $\dep(\Theta)\geqslant 2$ and consider the fixed $J\in I^\S$ such that $\dep(r_J(\Theta))<\dep(\Theta)$. In view of the results of section \ref{Case two or three}, $\p_J^\S(e_\Omega)$ is a linear combination of some elements $e_{\Omega'}$ where $\Omega'$ is an orbit included in $M_J(\Omega)$, which hence satisfies $\dep(\Omega') \geqslant \dep(\Omega) -p$, where $p = 1$, 2 or 3 when $J$ is of type A, B or C respectively. 

If $\dep(\Omega') \geqslant \dep(r_J(\Theta))$ and $\Omega'\neq r_J(\Theta)$ for every orbit $\Omega' \subseteq M_J(\Omega)$, we get by induction that $\p^\S_{g_{r_J(\Theta)}}(e_{\Omega'})$ is a linear combination of elements $e_{\Theta'}$ with $\Theta'$ distinct from $\Theta_{J_{r_J(\Theta)}} = \Theta_{J_\Theta}$, whence the result. The only obstructions to this situation are the two following configurations : 
\begin{itemize}
\item $M_J(\Theta)$ in configuration {\bf B3}, $M_J(\Omega)$ in configuration {\bf B4}, and $\dep(\Theta) = \dep(\Omega)$ maximal in $M_J(\Theta)$ and in $M_J(\Omega)$.
\item $M_J(\Theta) = M_J(\Omega)$ in configuration {\bf C5}, $\Omega$ of maximal depth, and $\Theta$ of second maximal depth, in $M_J(\Theta)$.
\end{itemize}

The second point occurs if $\G = D_4$ with $\Omega = \{\a_1+\a_2+\a_3+\a_4\}$ and $\Theta$, $J$ and $K$ as in the proof of (i) above (hence $g_\Theta = \D_K\D_J$), whence $\p^\S_{g_\Theta}(e_\Omega) = b^3\p^\S_K(e_{\Theta_K}) = 0$ and the result in that case. 

The first point occurs once if $\G = E_6$, for $J = \{3,5\}$, $\Omega = \{\a_2+\a_3+\a_4+\a_5\}$ and $\Theta = \{\a_1+\a_3+\a_4+\a_5, \a_3+\a_4+\a_5+\a_6\}$ with the notations of \cite[Planche V]{Bo}. But the possibilities for $(J_\Theta,g_\Theta)$ are then $(\{1,6\}, \D_J\D_{\{2\}}\D_{J})$, $(J, \D_{\{1,6\}}\D_{\{2\}}\D_{J})$ or $(J, \D_{\{2\}}\D_{\{1,6\}}\D_{J})$, whence $\p^\S_{g_\Theta}(e_\Omega) = b^3d^2e_{\{\a_{2}\}}$ and the result in that case.

The first point also occurs if $\G = A_{2n-1}$ ($n\geqslant 4$), for $J = \{n-j,n+j\}$, $\Omega = \{\sum_{i=n-j}^{n+j}\a_i\}$ and $\Theta = \{\sum_{i=n-3j}^{n-j}\a_i, \sum_{i=n+j}^{n+3j}\a_i\}$ where $1\leqslant j\leqslant \frac{n-1}{3}$, with the notations of \cite[Planche I]{Bo}. But in that case we necessarily get $g_\Theta = g_{r_Kr_J(\Theta)}\D_K\D_J$ where $K = \{n-k,n+k\}$ for some $j+1\leqslant k \leqslant 3j$. Hence the $K$-mesh $M_K(r_J(\Omega))$ is in configuration {\bf B2} and we get $\p_K^\S\p_J^\S(e_\Omega) = (bd)^2e_{r_J(\Omega)}$, with $\dep(r_J(\Omega)) = \dep(\Theta)-2 = \dep(r_Kr_J(\Theta))$ and $r_J(\Omega) \neq r_Kr_J(\Theta)$, whence the result in those cases by induction on $r_Kr_J(\Theta)$. 

Since there is no configuration {\bf B3} in $\Phi^+$ for $\G = D_n$, the proof is completed. \qed

\begin{thm}\label{thm irred}Assume that condition \emph{($\bigstar$)} holds : 
\begin{itemize}
\item[($\bigstar$)] $R$ is an integral domain, $f_i(e_{\a_i})\neq 0$ for all $i \in I$, and there exists a totally ordered integral domain $R_0$ and a ring homomorphism $\rho : R \to R_0$, such that $\rho(t) > 0$ for $t \in \{{a,b,c,d}\}$, and $\Image(f_i) \subseteq \ker(\rho)$ for all $i \in I$.
\end{itemize}
Then the LK-representations $\psi^\S$ and $\psi^\S_{\gr}$ are irreducible over $V_\K^\S$.
\end{thm}
\proof The result follows from propositions \ref{first part} and \ref{second part} : the second one applies with $Q = \ker(\rho)$, and the first one applies since under condition ($\bigstar$), we get $a$, $d$ and $f = f_i(e_{\a_i})$ non-zero in $R$, $f \in \ker(\rho)$ and $\rho(\check{d}) < 0$ (since $\rho(d\check{d}) = \rho(-bc) < 0$), so lemma \ref{tech result 1} shows that the coefficients $f_JQ_J(\p^\S_J)(e_{\Theta_K})$, for $m^\S_{J,K}\geqslant 3$, are non-zero in $R$ (look at the image of the parenthesized expressions by $\rho$). \qed

\begin{rem}This result extends to all the spherical and connected cases but $A_{2n}$ what was done for $\G = E_{6}$ in \cite[Section 6.1]{Ma}, by a totally different method. 

The question for $\G = A_{2n}$ ($n\geqslant 2$) is still open, as far as I know. The strategy used here does not directly apply because of the orbit $J$ of type D in $\G$ : first, lemma \ref{tech result 1} above is not true in this case (look at the orbit $\Theta'_J$), and second, the image of $P_J(\psi^\S_J)$ is not of dimension one, but two (at least when $f\neq 0$).

The case $\G = A_2$ is trivial however, since then $I^\S = \{I\}$, $V^\S_\K = \K{e_{\Theta_I}}\oplus \K{e_{\Theta'_I}}$, and $\psi^\S_I$ is the homothety of ratio $bcf$, hence in this case, $\psi^\S$ is not irreducible.
\end{rem}

\subsection{Non-equivalence}\label{Non-equivalence}\mbox{}\medskip

In corollary \ref{twisted simil} and theorem \ref{twisted equiv} below, we study the case of two twisted LK-representations $\psi^\S$ and $\psi'^\S$ of $B^+_{\G^\S} \approx (B_\G^+)^\S$ induced by two LK-representations $\psi$ and $\psi'$ of $B^+_\G$ for a \emph{fixed} Coxeter matrix $\G$.

On the contrary in proposition \ref{twisted equiv 2} below, we focus on the particular case of the Coxeter type $B_n$, $n\geqslant 3$, and study the case of three twisted LK-representations $\psi^\S$, $\psi'^{\S'}$ and $\psi''^{\S''}$ of $B^+_{B_n}$ induced by three LK-representations $\psi$, $\psi'$ and $\psi''$ of $B^+_{\G}$, $B^+_{\G'}$ and $B^+_{\G''}$ respectively, where $\G = A_{2n-1}$, $\G' = D_{n+1}$ and $\G'' = A_{2n}$.

\medskip

But let us begin with a result on the diagonalizability of the maps $\psi^\S_J$.

\begin{pro}\label{twisted diag}Let $J$ be an orbit of $I$ under $\S$. We set $\check{d} = a-d$.
\begin{enumerate}
\item If $J$ is of type A and if $d$, $\check{d}$ and $f$ are pairwise distinct, then $\psi^\S_J$ is diagonalizable over $\K$, with eigenvalues $d$, $\check{d}$ and $f$ of multiplicity $N_{\textrm{\bf A2}}+N_{\textrm{\bf A3}}$, $N_{\textrm{\bf A3}}$ and $N_{\textrm{\bf A1}}$ respectively.
\item If $J$ is of type B and if $d^2$, $d\check{d}$, $\check{d}^2$ and $df$ are pairwise distinct, then $\psi^\S_J$ is diagonalizable over $\K$, with eigenvalues $d^2$, $d\check{d}$, $\check{d}^2$ and $df$ of multiplicity $N_{\textrm{\bf B2}}+N_{\textrm{\bf B3}}+N_{\textrm{\bf B4}}$, $N_{\textrm{\bf B3}}+N_{\textrm{\bf B4}}$, $N_{\textrm{\bf B4}}$ and $N_{\textrm{\bf B1}}$ respectively.
\end{enumerate}
\end{pro}
\proof Let us prove (i). In view of proposition \ref{image des P(psiJ)}, the map $\psi^\S_J$ is diagonalizable over $\K$ and its eigenvalues can be $f_J(e_{\Theta_J}) = f$, $d$ or $\check{d}$. But we more precisely claim that for any $J$-mesh $M$ in configuration {\bf A1} (resp. {\bf A2}, resp. {\bf A3}), some suitable linear combinations of $e_{\Theta_J}$ and $e_\Theta$, $\Theta \subseteq M$, are eigenvectors of $\psi^\S_J$ for the value $f$ (resp. $d$, resp. $d$ and $\check{d}$), whence the result. 

The claim is obvious for {\bf A1} since $e_{\Theta_J}$ is clearly an eigenvector of $\psi^\S_J$ for the value $f$. For configuration {\bf A2}, the linear combination $e_\Theta + \lambda e_{\Theta_J}$ is suitable if $\lambda$ satisfies $f_J(e_\Theta) + \lambda f = d\lambda$, and such a $\lambda$ exists in $\K$ since $f \neq d$ by assumption. Now since $d\neq \check{d}$, the block of $\p^\S_J$ for a $J$-mesh $M$ in configuration {\bf A3} is certainly similar to the diagonal matrix $\textrm{Diag}(d,\check{d})$, \ie there certainly exists a linear combination $\tilde e$ of $e_\Theta$ and $e_{r_J(\Theta)}$, where $M = \Theta\cup r_J(\Theta)$, that is an eigenvector of $\p^\S_J$ for the eigenvalue $d$ (resp. $\check{d}$). The linear combination $\tilde e + \lambda e_{\Theta_J}$ is then an eigenvector of $\psi^\S_J$ for the eigenvalue $d$ (resp. $\check{d}$) whenever $f_J(\tilde e) + \lambda f = d\lambda$ (resp. $\check{d}\lambda$), and such a $\lambda$ exists in $\K$ since $f \neq d$ (resp. $f \neq \check{d}$) by assumption.

The proof of (ii) is similar (recall that there is no configuration {\bf B5} in $\Phi^+$).\qed

\begin{rem}One can prove an analogue of the previous results for the orbit $J$ of type C in $\G = D_4$ (when $|\S| = 3$), but we will not need it in what follows. The eigenvalues in this case would be $d^3$ (of multiplicity $N_{\textrm{\bf C2}}+N_{\textrm{\bf C5}} = 2$), and $d^2\check{d}$, $d\check{d}^2$, $\check{d}^3$ and $d^2f$ (of multiplicity $N_{\textrm{\bf C5}} = N_{\textrm{\bf C1}} = 1$).

There is also an analogue for the orbit $J$ of type D in $\G = A_{2n}$, if we assume that $-d\check{d} = bc$ is a square in $\K$, so that the factor $X^2+(d\check{d})^3$ of $P_J$ splits in $\K[X]$. 
\end{rem}

Now let us fix another LK-representation $\psi' : B^+ \to \EL(V)$ associated with a quadruple $(a',b',c',d')\in R^4$ (such that $d'^2-a'd'-b'c' = 0$) and an LK-family $(f'_{i})_{i\in I}$ relatively to $(a',b',c',d')$. We assume that $b'$, $c'$ and $d'$ are non-zero in $R$.

As in the preamble of the current section, the LK-family $(f'_{i})_{i\in I}$, seen as an element of $(V_\K^\star)^I$, is entirely determined by the common value $f'$ of the $f'_{i}(e_{\a_i})$ for $i\in I$, and the LK-representation $\psi'$ induces a twisted LK-representation $\psi'^\S : (B^+)^\S \to \EL(V_\K^\S)$ of $(B^+)^\S$.

\begin{cor}\label{twisted simil}We set $\check{d} = a-d$ and $\check{d'} =  a'-d'$. 
\begin{enumerate}
\item Assume that $\G \neq A_{2n}$, $A_3$ and that $d$, $\check{d}$ and $f$ (resp. $d'$, $\check{d'}$ and $f'$) are pairwise distinct. Then for any orbit $J \in I^\S$ of type A, the maps $\psi^\S_J$ and $\psi'^\S_J$ are similar over $\K$ if and only if $(d,\check{d},f) = (d',\check{d'},f')$.
\item Assume that $\G = A_{2n}$, $n \geqslant 3$, and that $d^2$, $d\check{d}$, $\check{d}^2$ and $df$ (resp. $d'^2$, $d'\check{d'}$, $\check{d'}{}^2$ and $d'f'$) are pairwise distinct. Then for any orbit $J$ of type B, the maps $\psi^\S_J$ and $\psi'^\S_J$ are similar over $\K$ if and only if $(d,\check{d},f) = \pm(d',\check{d'},f')$.
\end{enumerate}

In particular, if $\G \neq A_{2n}$, $A_3$ and $(d,\check{d},f) \neq (d',\check{d'},f')$ (resp. if $\G = A_{2n}$, $n\geqslant 3$, and $(d,\check{d},f) \neq \pm(d',\check{d'},f')$), then the twisted LK-representations $\psi^\S$ and $\psi'^\S$ are not equivalent over $\K$.
\end{cor}
\proof Of course if $\psi^\S$ and $\psi'^\S$ are equivalent, then the maps $\psi^\S_J$ and $\psi'^\S_J$ are similar for every $J\in I^\S$. But in view of proposition \ref{twisted diag} above, the maps $\psi^\S_J$ and $\psi'^\S_J$ are similar if and only if they have the same eigenvalues, with same multiplicity. 

In case (i), the three multiplicities are $N_{\textrm{\bf A2}}+N_{\textrm{\bf A3}} > N_{\textrm{\bf A3}} > N_{\textrm{\bf A1}}$ (the second inequality is strict since we avoid the case $\G = A_3$), whence the result.

In case (ii), the four multiplicities are $N_{\textrm{\bf B2}}+N_{\textrm{\bf B3}}+N_{\textrm{\bf B4}} > N_{\textrm{\bf B3}}+N_{\textrm{\bf B4}} > N_{\textrm{\bf B4}} = N_{\textrm{\bf B1}}$ (the first inequality is strict since we avoid the case $\G = A_4$), so $\psi^\S_J$ and $\psi'^\S_J$ are similar if and only if $d^2 = d'^2$, $d\check{d} = d'\check{d'}$ and $\{\check{d}^2,df\} = \{\check{d'}{}^2,d'f'\}$. And these three equalities are clearly equivalent to $(d,\check{d},f) = \pm(d',\check{d'},f')$. \qed

\begin{rem}\label{twisted simil 2}The case $\G = A_3$ (resp. $A_4$) can be dealt with as in case (i) (resp. (ii)) of the previous proposition and its proof, but the criterion of similarity of $\psi^\S_J$ and $\psi'^\S_J$ for $J$ of type A (resp. B) as then to be relaxed to $d = d'$ and $\{\check{d},f\} = \{\check{d'},f'\}$ (resp. $\{d^2,d\check{d}\} = \{d'^2,d'\check{d'}\}$ and $\{\check{d}^2,df\} = \{\check{d'}{}^2,d'f'\}$).

The case $\G = A_2$ is trivial since then $I^\S = \{I\}$, $V^\S_\K = \K{e_{\Theta_I}}\oplus \K{e_{\Theta'_I}}$, and $\psi^\S_I$ (resp. $\psi'^\S_I$) is the homothety of ratio $bcf$ (resp. $b'c'f'$). Hence in this case, the twisted LK-representations $\psi^\S$ and $\psi'^\S$ are equivalent if and only if $bcf = b'c'f'$.
\end{rem}

\begin{thm}\label{twisted equiv}Assume that $\G \neq A_2$ and that condition \emph{($\bigstar$)} holds for both $\psi$ and $\psi'$, for the \emph{same} morphism $\rho : R \to R_0$ if $\G = A_3$ or $\G = A_{2n}$. Then if $(d,\check{d},f) \neq (d',\check{d'},f')$, the twisted LK-representations $\psi^\S$ and $\psi'^\S$ are not equivalent over $\K$.
\end{thm}
\proof Under condition ($\bigstar$), the elements $d,\, \check{d}$ and $f$ are pairwise distinct since their image by $\rho$ are positive, negative and zero respectively (for $\check{d}$, the result follows from the identity $d\check{d} = -bc$). Similarly, the elements $df$, $d^2$, $d\check{d}$ and $\check{d}^2$ are pairwise distinct : indeed, we get $\rho(df) = 0$, $\rho(d\check{d}) < 0$, $\rho(d^2)>0$ and $\rho(\check{d}^2) > 0$, and $d^2 \neq \check{d}^2$ since $d^2 - \check{d}^2 = (d-\check{d})a$ and hence $\rho(d^2 - \check{d}^2)>0$. 

The analogue occurs for $\psi'$ and hence the results of corollary \ref{twisted simil} and remark \ref{twisted simil 2} applies. This gives the result if $\G \neq A_{2n}$, $A_3$. If $\G = A_{2n}$, $n\geqslant 3$, then the condition stated in corollary \ref{twisted simil} is $(d,\check{d},f) \neq \pm(d',\check{d'},f')$, but we clearly have $(d,\check{d},f) \neq -(d',\check{d'},f')$ since $\rho(d)$ and $\rho(d')$ are positive, whence the result in that case. Similarly for $\G = A_3$ (resp. $A_4$), the condition of non-equivalence derived from remark \ref{twisted simil 2} should be $d \neq d'$ or $\{\check{d},f\} \neq \{\check{d'},f'\}$ (resp. $\{d^2,d\check{d}\} \neq \{d'^2,d'\check{d'}\}$ or $\{\check{d}^2,df\} \neq \{\check{d'}{}^2,d'f'\}$). But here $\{\check{d},f\} = \{\check{d'},f'\}$ is equivalent to $(\check{d},f) = (\check{d'},f')$ since $\rho(d)$ and $\rho(d')$ are positive and $\rho(f) = \rho(f') = 0$, whence the result for $\G = A_3$. Finally for $\G = A_4$, since the images by $\rho$ of $d^2$, $\check{d}^2$, $d'^2$ and $\check{d'}{}^2$ (resp. of $d\check{d}$ and $d'\check{d'}$, resp. of $df$ and $d'f'$) are positive (resp. negative, resp. zero), the condition $\{d^2,d\check{d}\} = \{d'^2,d'\check{d'}\}$ and $\{\check{d}^2,df\} = \{\check{d'}{}^2,d'f'\}$ of remark \ref{twisted simil 2} is equivalent to $(d^2,d\check{d},\check{d}^2,df) = (d'^2,d'\check{d'},\check{d'}{}^2,d'f')$, \ie $(d,\check{d},f) = (d',\check{d'},f')$ since, again, $\rho(d)$ and $\rho(d')$ are positive. \qed

\medskip

Let us finally turn to the case of the Coxeter type $B_n$ by considering three twisted LK-representations $\psi^\S$, $\psi'^{\S'}$ and $\psi''^{\S''}$ of $B^+_{B_n}$ induced by three LK-representations $\psi$, $\psi'$ and $\psi''$ of $B^+_{\G}$, $B^+_{\G'}$ and $B^+_{\G''}$ respectively, where $\G = A_{2n-1}$, $\G' = D_{n+1}$ and $\G'' = A_{2n}$.

\medskip

The twisted LK-representation $\psi''^{\S''}$ is clearly non-equivalent to the two others since it is of dimension $|\Phi^+_{\G''}/\S''| = n(n+1)$ whereas the two others are of dimension $|\Phi^+_{\G}/\S| = |\Phi^+_{\G'}/\S'| = n^2$. The following proposition proves that $\psi^\S$ ans $\psi'^{\S'}$ are not equivalent, under the usual assumptions on the parameters $(a,b,c,d)$ and $(f_i)_{i\in I}$ of $\psi$, and $(a',b',c',d')$ and $(f'_i)_{i\in I'}$ of $\psi'$. We set as above $\check{d} = a-d$ and $f = f_i(e_{\a_i})$ for any $i\in I$.

\begin{pro}\label{twisted equiv 2}Assume that $d,\, \check{d}$ and $f$ (resp. $df$, $d^2$, $d\check{d}$ and $\check{d}{}^2$) are pairwise distinct, and that the analogue holds for $\psi'$, which is the case for example if condition \emph{($\bigstar$)} holds for both $\psi$ and $\psi'$. Then the twisted LK-representations $\psi^\S$ and $\psi'^{\S'}$ are not equivalent over $\K$.
\end{pro}
\proof The fact that condition ($\bigstar$) implies the given conditions on the parameters has been shown in the proof of theorem \ref{twisted equiv 2}. But then proposition \ref{twisted diag} applies and shows that the images by $\psi$ and $\psi'$ of a given standard generator of $B^+_{B_n}$ do not have the same number of eigenvalues, whence the result. \qed

\bigskip

\begin{flushleft}{\small 
\textsc{Universit\'e Paris-Sud 11, Math\'ematiques - B\^at. 425, 91\,405 Orsay, France.}\\
\emph{E-mail address} : \url{anatole.castella@math.u-psud.fr}
}
\end{flushleft}

\end{document}